\documentclass[12pt,leqno]{article}
\usepackage{latexsym}
\usepackage[all]{xy}
\usepackage{amsmath}
\usepackage{amssymb}
\usepackage{amsfonts}
\usepackage{color}
\usepackage{theorem}

\usepackage{hyperref}
\usepackage{xcolor}
\hypersetup{colorlinks=true,urlcolor=blue,linkcolor=blue,citecolor=blue}

\newcommand{\op}[1]{\operatorname{#1}}

\newcommand {\Perf} {\mathsf{Perf}}
\newcommand{\iPerf}{\mathit{Perf}}
\newcommand {\Vect} {\mathsf{Vect}}

\newcommand {\Map} {\mathbf{Map}}
\newcommand {\rh} {\mathbb{R}\underline{Hom}}

\newcommand {\OO} {\mathcal{O}}
\newcommand {\DR}{\mathbf{DR}}

\newcommand {\fb} {\widehat{\partial}}
\newcommand {\A} {\mathcal{A}}
\newcommand {\hH} {\mathcal{H}}
\newcommand {\T} {\mathbb{T}}

\newcommand  {\dg}     {\mathbf{dg}}
\newcommand  {\dgcat}     {\mathbf{dgCat}}

\newcommand  {\gdg}     {\mathbf{dg}^{gr}}

\newcommand  {\medg}     {\mathbf{dg}_{\epsilon}^{gr}}

\newcommand  {\dgart}     {\mathbf{dgArt}}

\newcommand  {\dAff}     {\mathbf{dAff}}

\newcommand  {\dSt}   {\mathbf{dSt}}
\newcommand  {\LL}   {\mathcal{L}}

\newcommand{\s}{\infty}
\newcommand{\HH}{\mathbb{H}}
\newcommand{\HHc}{\widetilde{\HH}_c}
\newcommand{\HHcd}{\widetilde{\HH}_{c,DR}}

\newcommand{\D}{\mathcal{D}}

\newcommand{\Xbar}{\mathfrak{X}}
\newcommand{\Dqcoh}{\mathsf{D}_{\text{qcoh}}}
\newcommand{\bepsilon}{\boldsymbol{\epsilon}}
\newcommand{\barE}{\mathcal{E}}
\newcommand{\jfrak}{\mathfrak{j}}
\newcommand{\barV}{{\mathcal{V}}}
\newcommand{\barnabla}{\mathfrak{d}}
\newcommand{\barW}{{\mathcal{W}}}
\newcommand{\Jfb}{\mathfrak{J}}

\newtheorem{introthm}{Theorem}
\newtheorem{thm}{Theorem}[section]
\newtheorem{prop}[thm]{Proposition}
\newtheorem{lem}[thm]{Lemma}
\newtheorem{sublem}[thm]{Sublemma}
\newtheorem{df}[thm]{Definition}
\newtheorem{cor}[thm]{Corollary}

{\theorembodyfont{\rmfamily} \newtheorem{rmk}[thm]{Remark}}

\numberwithin{equation}{section}

\newcommand{\Appendix}[1]{%
  \refstepcounter{section}%
  \addtocontents{toc}{\protect\setcounter{tocdepth}{1}}
  \addcontentsline{toc}{section}%
    {\bfseries\appendixname~\thesection:\ #1}%
    {\medskip\noindent \Large\bfseries\appendixname\ \thesection:\ #1}%
\sectionmark{#1}\smallskip\noindent
\renewcommand{\theequation}{{\bf 
{{\thesection}}.{\arabic{equation}}}}
}

\begin{document}

\title{\textbf{Moduli of flat connections on smooth varieties}}  

\author{Tony Pantev and Bertrand To\"{e}n}


\date{September  2021}

\maketitle

\begin{abstract} This paper is a companion to \cite{ptbetti}. We study
  the moduli functor of flat bundles on smooth, possibly non-proper,
  algebraic variety $X$ over a field of characteristic zero. For
  this we introduce the notion of a \emph{formal boundary} of $X$,
  denoted by $\fb X$, which is a formal analogue of the boundary at
  $\s$ of the Betti topological space associated to $X$. We explain
  how to construct two derived moduli functors $\Vect^{\nabla}(X)$ and
  $\Vect^{\nabla}(\fb X)$, of flat bundles on $X$ and on $\fb X$, as
  well as a restriction map
  $\mathsf{R} : \Vect^{\nabla}(X) \longrightarrow \Vect^{\nabla}(\fb X)$.

  This work contains two main results. First we prove that the
  morphism $\mathsf{R}$ comes equipped with a canonical shifted
  Lagrangian structure. This result can be understood as the de Rham
  analogue of the existence of Poisson structures on moduli of local
  systems on $X$. As a second statement, we prove that the geometric
  fibers of $\mathsf{R}$ are representable by \emph{quasi-algebraic
    spaces}, a slight weakening of the notion of algebraic spaces.
\end{abstract}


\tableofcontents

\section*{Introduction}
\addcontentsline{toc}{section}{Introduction}

This work is a sequel of \cite{ptbetti} in which we studied moduli of
local systems on a topological space underlying a smooth non-proper
complex algebraic variety $X$. One of the main results of
\cite{ptbetti} asserts that this moduli is a derived Artin stack
endowed with a natural shifted Poisson structure whose symplectic
leaves can be studied by fixing monodromies of local systems at
infinity.

In this paper we begin the study of the de Rham analogue of the
results of \cite{ptbetti}.  The content of the present work can be
summarized in the statement that for a smooth variety $X$ over a field
$k$ of characteristic $0$ the deived moduli $\Vect^{\nabla}(X)$ of
flat connections on $X$ carries a canonical shifted Poisson
structure. However, this statement needs to be qualified as
$\Vect^{\nabla}(X)$ is not representable for non-proper $X$. Thus in
order to state and prove the existence of Poison structures in this
context we have to overcome a number of technical difficulties.

When restricted to regular connections, it is possible
to approach this representability question by working on some good
compactification, as done in \cite{ni}. What we propose
in this paper is slightly different, as we do not assume
any regularity assumptions and also propose an intrinsic construction, 
independent of any choice of compactification.
A key ingredient for this  is the notion of a \emph{\bfseries
  formal boundary} $\fb X$ of a smooth variety $X$. It is difficult to
make sense of the formal boundary directly as a geometric object but
one can make sense of it as a non-commutative space. In particular it
is possible to define $\s$-categories of vector bundles and flat
bundles on $\fb X$. The putative object $\fb X$ is morally the
punctured formal completion of $\Xbar$ along $D$, for $\Xbar$ a smooth
compactification of $X$ with $D = \Xbar - X$ a normal crossings
divisor.  Rigid analytic and formal versions of $\fb X$ have been
considered previously in \cite{bete,ef,hpv}. The novelty here is the
systematic study of corresponding  de Rham theory: vector bundles with
connections on $\fb X$ and their de Rham complexes. Some glimpses of
such theory already exist the literature. In particular, when $X$ is a
curve, the de Rham theory of $\fb X$ was developed and analyzed by
S.~Raskin in \cite{ra} in the context of the local geometric Langlands
correspondence. In this paper we deal with the case of a higher
dimensional $X$. We construct derived stacks $\Vect^{\nabla}(X)$ and
$\Vect^{\nabla}(\fb X)$ of flat bundles on $X$ and $\fb X$, together
with a restriction map
$\mathsf{R} : \Vect^{\nabla}(X) \longrightarrow \Vect^{\nabla}(\fb
X)$. We study the infinitesimal properties of these derived stacks,
and show in particular that they are formally representable at any
field valued point. The formal representability allows us to define
the notion of shifted symplectic and shifted Lagrangian structures on
these derived stacks, even though they are not representable. Our
first main result then is the following theorem.

\

\begin{introthm}\label{ti1}
There exists a canonical $(3-2d)$-shifted Lagrangian structure 
on the restriction map
$$
\mathsf{R} : \Vect^{\nabla}(X) \longrightarrow \Vect^{\nabla}(\fb X).
$$
\end{introthm} 

\

\medskip

\noindent
At the linear level of tangent complexes, the above theorem is an
incarnation of Poincar\'e duality in de Rham cohomology and de Rham
cohomology with compact supports.  The existence of the Lagrangian
structure globally is itself a version of Poincar\'e duality relative
to various derived base schemes, together with the general existence
result \cite[Theorem~3.7]{to3}.  Also, Theorem \ref{ti1} \ immediately implies
the existence of a $(2-2d)$-shifted Poisson structure on
$\Vect^{\nabla}(X)$ thanks to the comparison result 
\cite[Theorem~4.22]{mesa}.

Our second main result in this work is the following representability
result.  We fix a flat bundle at infinity
$V_\s \in \Vect^{\nabla}(\fb X)(k)$ and consider the fiber of
$\mathsf{R}$ at $V_\s$ denoted by $\Vect ^{\nabla}_{V_\s}(X)$. Our
original goal was to prove that $\Vect ^{\nabla}_{V_\s}(X)$ is
representable by a derived Artin stack (even algebraic space if no
components of $X$ are proper) locally of finite presentation over
$k$. Though we have not been able to prove this last statement, we
prove the following weaker version.

\

\begin{introthm}\label{ti2}
The derived stack $\Vect^\nabla_{V_\s}(X)$
is a derived quasi-algebraic space locally of finite presentation 
in the sense of definition \ref{d17}.
\end{introthm}

\

\medskip

\noindent
Derived quasi-algebraic spaces are almost algebraic spaces - they
satisfy all conditions in the Artin-Lurie representability criterion
(see \cite{lu2} and Appendix~\ref{sec:artin-lurie}) except that they
may not be locally of finite presentation as a
functor. Quasi-algebraic spaces only satisfy local finite
presentability generically, and the result is that these derived
stacks only have a smooth atlas generically, i.e. have a smooth atlas
whose image is Zariski dense in an appropriate sense.

\

\medskip

\noindent
{\bfseries Acknowledgements:} We would like to thank Sasha Efimov,
Dmitry Kaledin, and Gabriele Vezzosi for several illuminating
discussions on the subject of this work.

During the preparation of this work Bertrand To\"en was partially
supported by ERC-2016-ADG-741501.  Tony Pantev was partially supported
by NSF research grants DMS-1601438 and DMS-1901876, by Simons
Collaboration grant \# 347070, and by the Laboratory of Mirror
Symmetry NRU HSE, RF Government grant, ag. No 14.641.31.0001.

\

\subsection*{Notations and conventions}
\addcontentsline{toc}{subsection}{Notation and conventions}

\begin{description}
\item[$k$] - a field of characteristic zero.
\item[$G$] - an affine reductive group over $k$
\item[$X$] - a smooth variety over $k$.
\item[$\fb X$] - the formal boundary of $X$.
\item[$\Vect^{\nabla}(X)$] - the derived stack of flat connections on $X$.
\item[$\Vect^{\nabla}(\fb X)$] - the derived stack of flat connections
  on $\fb X$.
\item[$\mathsf{R} : \Vect^{\nabla}(X) \longrightarrow
  \Vect^{\nabla}(\fb X)$] - the restriction to the boundary map.
\item[$\Xbar$] - a smooth compactification of $X$ with a normal
  crossings divisor boundary $D = \Xbar - X$.
\item[$R_{G}(\pi_{1}(X,x))$] - the $G$-character scheme parametrizing
  representations of $\pi_{1}(X,x)$ in $G$.
\item[$Loc_{G}(X)$] - the derived moduli stack of $G$-local systems on $X$.
\item[$\T$] - the $\s$-category of spaces, or equivalently the
  $\s$-category of simplicial sets.
\item[$t_{0}Loc_{G}(X)$] - the underived truncation of $Loc_{G}(X)$. 
\item[$\mathsf{commalg}_{k}$] - the category of commutative
  $k$-algebras.
\item[$\mathsf{C}_{\lambda}$] - the conjugacy class of a group element
  $\lambda \in G$,
\item[$Loc_{G}(X,\lambda)$] - the derived moduli stack of $G$ local
  systems on $X$ with monodromy at infinity in $\mathsf{C}_{\lambda}$
  (assumes $X$ admits a smooth compactification with a smooth
  connected divisior at infinity).
\item[$G*G$] - the derived commuting variety of $G$.
\item[$\D_{X}$] - the sheaf of rings of differential
  operators on a smooth variety $X$.
\item[$B$] - a connective cdga.
\item[$\Dqcoh(\D_{X,B})$] - the dg-category of all $\D_{X}\otimes_{k}
  B$-modules, whose underlying $\mathcal{O}_{X}\otimes_{k} B$-modules
  are quasi-coherent on $X\times Spec\ B$.
\item[\emph{a symmetric monoidal dg-category}] - 
an $E_\s$-algebra object inside the symmetric monoidal $\s$-category
of locally presentable dg-categories (see \cite[Section~2]{azto}).
\item[$k-\dg$] - the $\s$-category of complexes of $k$-modules.
\item[$k-\medg$] - the $\s$-category of graded mixed $k$-modules.
\item[$\DR_X-\medg$] - dg-category of cofibfrant graded mixed
  $\DR_{X}$-dg-modules for a smooth affine variety $X$ over $k$.
\item{$\hH$} - the group stack of autoequivalences of
  $B\mathbb{G}_{a}$.
\item[$\mathbb{D}g^{lp}$] - the derived stack of locally presentable
  dg-categories.
\item[$\mathbb{D}g^{lp}(B\hH)$] - the $\s$-category of
  $\hH$-equivariant locally presentable dg-categories.
\item[$\jfrak_n : \Xbar_{(n)}:= Spec\ (\OO_{\Xbar}/I_D^n)
\to \Xbar$] -
the $(n-1)$-th infinitesimal thickening of $D$ inside $\Xbar$
\item[$\widehat{\jfrak} : \widehat{\Xbar} \to \Xbar$] - the full formal
  neighborhood of $D$ inside $\Xbar$.
\item[$\Perf(\widehat{\Xbar})$] - the derived stack of perfect
  complexes on $\widehat{\Xbar}$.
\item[$\Perf(\fb X)$] - the derived stack of perfect
  complexes on the formal boundary of $X$.
\item[$\Perf^{ex}(\fb X)$] - the derived stack of extendable perfect
  complexes on the formal boundary of $X$.
\item[$\Perf^{\nabla}(\widehat{\Xbar})$] - the derived stack of perfect
  complexes of flat bundles on $\widehat{\Xbar}$.
\item[$\Perf^{\nabla}(\fb X)$] - the derived stack of perfect
  complexes of flat bundles on the formal boundary of $X$.
\item[$\Perf^{\nabla,ex}(\fb X)$] - the derived stack of extendable perfect
  complexes of flat bundles on the formal boundary of $X$.
\end{description}

\section{A brief review of the Betti case}

Before we plunge into the technical aspects of formal boundaries and
their de Rham theory, it is useful to review the Betti case results from
\cite{ptbetti} which motivate the present discussion.

A precursor of this whole line of investigation is the classical story
about Poisson structures on the moduli of representations of
fundamental groups of topological surfaces: if $X$ is a compact
oriented topological surface, and $G$ is a complex reductive group,
then it is well known (see \cite{fock-rosly,
  ghjw,goldman,guruprasad-rajan}) that the smooth part of the moduli
space of representations $\rho : \pi_{1}(X) \to G$ carries a canonical
algebraic Poisson structure, and that the the symplectic leaves of
this Poisson structure are moduli spaces of representations $\rho$
whose values at the loops at infinity belong to fixed conjugacy
classes in $G$. In \cite{ptbetti} we extended this story to higher
dimensional smooth open varieties $X$ simultaneously refining the
moduli problem for representations of the fundamental group in a way
which removes the smoothness restriction on the moduli. More precisely
we proved the following

\

\begin{thm} {\em \cite[Theorem~4.7]{ptbetti}}  \label{thm-betti}
  Fix a field $k$ of
  $\mathsf{char} k = 0$.  Let $X$ be a $d$-dimensional smooth complex
  algebraic variety and let $G$ be a reductive algebraic group over
  $k$. Then
\begin{itemize}
\item[\textcolor{blue}{{\bfseries (1)}}] The derived moduli stack
  $Loc_{G}(X)$ of $G$-local systems
  on $X$ has a natural $(2-2d)$-shifted Poisson structure.
\item[\textcolor{blue}{{\bfseries (2)}}] This shifted Poisson
  structure admits generalized symplectic leaves.  Among those are the
  derived moduli of $G$ local systems with fixed monodromy at
  infinity.
\end{itemize}
\end{thm}

\

\begin{rmk}
\begin{enumerate}
\item[(i)] When $d = 1$ the Poisson structure in
  Theorem~\ref{thm-betti} (1)  specializes to Goldman's Poisson
  structure on the moduli of representations $\pi_{1}(X) \to G$.
\item[(ii)] The precise formulation of the statement
  Theorem~\ref{thm-betti} (2) is tricky since for this one needs to
  understand how to fix local monodromies in the derived setting. To
  solve this problem in \cite{ptbetti} we had to deal with a couple of
  subtle issues:
\begin{itemize}
\item the fixing of the local monodromies can not be seen solely on
  underlying underived moduli stack $t_{0} Loc_{G}(X)$ and involves
  higher homotopy coherences;
\item in higher dimension an additional constraint called \emph{\bfseries
        strictness} has to be imposed on the local monodromies at
      infinity in order to select a symplectic leaf.
    \end{itemize}
\end{enumerate}
\end{rmk}

\

\noindent
In more details, suppose $X$ is a finite CW complex and $G$ a
reductive group over $k$. The derived moduli stack $Loc_{G}(X)$ of
$G$-local systems (=locally constant principal $G$ bundles) is a
derived Artin stack locally of finite presentation over $k$. The
truncated underived stack $t_{0} Loc_{G}(X)$ depends only on the
fundamental group of $X$. It is the moduli stack of representations of
$\pi_{1}(X,x)$ into $G$, i.e.
  \[
  t_{0} Loc_{G}(X) = \left[\left.
    R_{G}(\pi_{1}(X,x))\right/ G\right]
 \]
where $R_{G}(\pi_{1}(X,x))$ is the
 \emph{\bfseries character scheme} of $X$, namely
 the affine $k$-scheme representing the functor
 \[
 \xymatrix@R-2pc{
   R_{G}(\pi_{1}(X,x)) :
   \hspace{-1.5pc} & \mathsf{commalg}_{k} \ar[r] & \mathsf{Sets}, \\
& A \ar[r] & \op{Hom}_{\text{grp}}\left(\pi_{1}(X,x),G(A)\right).}
 \]
In this sense $Loc_{G}(X)$ refines the moduli of representations of
$\pi_{1}(X,x)$. Note however (see \cite[Section~1.2]{ptbetti}) that in
general the derived structure on $Loc_{G}(X)$ depends on the full
homotopy type of $X$ so it captures more information, than the moduli
of representations.

To explain the content of Theorem~\ref{thm-betti} properly, recall
that the \emph{\bfseries boundary of a topological space} $X$ is the
pro-homotopy type $\partial X:= \underset{K \subset X}{\lim}(X-K)$,
where the limit is taken in the $\s$-category $\T$ of spaces and over
the opposite category of compact subsets $K \subset X$. In general
$\partial X$ can be quite complicated but if $X$ is a the underlying
topological space of a smooth $d$-dimensional complex algebraic
variety, then $\partial X$ is equivalent to a constant pro-object in
$\T$ which has the homotopy type of a compact oriented topological
manifold of dimension $2d-1$. This implies (see \cite{ca} and
\cite[Section~4.1]{ptbetti}) that 
the canonical map
$\partial X \longrightarrow X$ induces a restriction morphism of
derived locally f.p. Artin stacks
$$
r : Loc_{G}(X) \longrightarrow Loc_{G}(\partial X).
$$
which is equipped with a canonical $(2-2d)$-shifted Lagrangian
structure with respect to the canonical shifted symplectic structure
on $Loc_{G}(\partial X)$.  Therefore by \cite[Theorem~4.22]{mesa} the map $r$ can be
viewed as a $(2-2d)$-shifted Poisson structure on $Loc_{G}(X)$, which
gives part (1) of Theorem~\ref{thm-betti}.

For part (2) of Theorem~\ref{thm-betti} the restriction on the local
monodromies at infinity has to be made precise in at least two
different ways. First we need to take into account the fact that the
boundary of $X$ has higher dimension and the local monodromy around
the boundary may get twisted as we move along a connected component of
the boundary. Second we have to make sure that when we have good
control over the interaction of the local monodromies going around
intersecting divisor components in a normal crossings
compactification.

To sketch how one deals with the first issue suppose that we have a smooth
compactification $\Xbar$ of $X$. To simplify the discussion assume
that $D = \Xbar - X$ is a smooth connected divisor. Then $\partial X$
has the homotopy type of an oriented circle bundle over $D$ classified
by $\alpha = c_{1}(N_{D/\Xbar}) \in H^{2}(D,\mathbb{Z})$. Now given
$\lambda \in G$ with centralizer $Z < G$, the group $S^{1}$ acts on
$BZ$ (via $\lambda$) and naturally on the conjugation quotient $[G/G]$
so that the $1$-shifted Lagrangian structure on the map $BZ \to [G/G]$ is
$S^{1}$-equivariant. Twisting by $\alpha$ gives a $1$-shifted
Lagrangian morphism
\begin{equation} \label{eq-twisted} 
  {}_{\alpha}\widetilde{BZ} \longrightarrow
  {}_{\alpha}\widetilde{[G/G]}
\end{equation}
of  locally constant families of derived Artin
stacks over $D$. Passing to global sections gives moduli stacks
\[
\begin{aligned}
Loc_{G}(\partial X) &  = \mathsf{Map}\left(\partial X,
  BG\right) =
  \Gamma\left(D_{i},{}_{\alpha}\widetilde{[G/G]}\right); \\
Loc_{Z,\alpha}(D) & =
\Gamma\left(D,{}_{\alpha}\widetilde{BZ}\right)
\end{aligned}
\]
of $G$ local systems on  $\partial X$ and of $\alpha$-twisted 
$Z$-local systems on $D$ respectively. Furthermore, since $D$
is a compact topological manifold endowed
  with a canonical orientation the map \eqref{eq-twisted} induces a
  $(3-2d)$-shifted Lagrangian morphism of derived Artin stacks
$$
Loc_{Z,\alpha}(D) \longrightarrow Loc_{G}(\partial X).
$$
By the Lagrangian intersection theorem of \cite[Section~2.9]{ptvv} the fiber
product of derived stacks
$$
Loc_{G}(X,\lambda) := Loc_{Z,\alpha}(D)
\underset{Loc_{G}(\partial X)}{\times}  Loc_{G}(X)
$$
has a canonical $(2-2d)$-shifted symplectic structure. By construction
\begin{itemize}
\item $Loc_{G}(X,\lambda)$ is the derived
stack of $G$-local systems  on $X$ whose local monodromy
around $D$ is fixed to be in the conjugacy class
$\mathsf{C}_{\lambda}$ of $\lambda$.
\item The  natural map
$$
Loc_{G}(X,\lambda) \longrightarrow Loc_{G}(X)
$$
thus realizes 
$Loc_{G}(X,\lambda)$
as a generalized  symplectic leaf of the $(2-2d)$-shifted
Poisson structure on $Loc_{G}(X)$. 
\end{itemize}    
This explains  Theorem~\ref{thm-betti} (2) in the case when $X$ admits a
compactification with a smooth divisor boundary.

To sketch how one deals with the second issue start again with a
smooth compactification $\Xbar$ of $X$ but this time assume that the
divisor at infinity $D = \Xbar - X = D_{1}\cup D_{2}$ has two smooth
irreducible components meeting transversally at a smooth connected
subvariety $D_{12}$. In this case
\[
\partial X \simeq \partial_{1} X \sqcup_{\partial_{12}X}\partial_{2} X.
\]
where $\partial_{i} X$ is an oriented circle bundle over
$D_{i}^{o}=D_{i}-D_{12}$, and $\partial_{12}X$ is an oriented $S^1
\times S^1$-bundle over $D_{12}$. Note that here each $\partial_{i} X$
has the homotopy type of an oriented compact manifold of dimension
$2d-1$ with boundary canonically equivalent to $\partial_{12}X$.

The problem we need to solve now is to understand what conditions on a
pair of commuting elements $\lambda_{1}, \lambda_{2} \in G$ will
guarantee that prescribing the $\lambda_{i}$ as the local monodromies
around the components $D_{i}$ will select a symplectic leaf in
$Loc_{G}(X)$. In \cite{ptbetti} we found a natural sufficient
condition called \emph{\bfseries strictness}. 

\

\begin{df} \label{df-strict}
A pair of commuting
  elements $(\lambda_1,\lambda_2) \in G \times G$ is called
  \emph{\bfseries strict} if the morphism
$$
BZ_{12} \longrightarrow [Z_1/Z_1]\times_{[G*G/G]}[Z_2/Z_2]
$$
is Lagrangian (for its canonical isotropic structure).
\end{df}

\

\noindent
Here $Z_{i}$ denotes the centralizer subgroup of $\lambda_{i}$,
$Z_{12}$ denotes the centralizer of the pair
$(\lambda_{1},\lambda_{2})$, \linebreak $G*G \subset G\times G$ is the
derived commuting variety of $G$, and $[G*G/G]$ is the stack quotient of $G*G$
by the diagonal conjugation action of $G$.

\

\begin{rmk} \label{rmk-strictness}
By definition strictness is a group theoretic property and in fact can
be expressed in elementary group theoretic terms. In
\cite[Proposition~4.9]{ptbetti} we show that if
$(\lambda_1,\lambda_2)$ is a commuting pair of elements in $G$, and
$u:=\mathsf{Id}-\mathsf{ad}(\lambda_1)$ and
$v:=\mathsf{Id}-\mathsf{ad}(\lambda_2)$ are the corresponding
endormorphisms of $\mathfrak{g}$, then the pair
$(\lambda_1,\lambda_2)$ is strict if and only $u$ is strict with
respect to the kernel of $v$, i.e. if and only if
$\op{Im}(v_{|\ker(u)})=\op{Im}(v)\cap \ker(u)$.
\end{rmk}

\

\noindent
With the notion of strictness at hand we can formulate a precise
version of the statement of Theorem~\ref{thm-betti} (2) in the case of
a strict normal crossings boundary divisor with two components:

\

\begin{thm} \label{thm-strict} {\em \cite[Theorem~4.7]{ptbetti}}
Let  $(\lambda_{1},\lambda_{2})$  be a 
strict pair of commuting elements in $G$ and let  
$$
Loc_{G}(X,\{\lambda_1,\lambda_2\})
$$
be the derived Artin stack of local systems on $X$ whose local
monodromy around $D_{i}$ belongs to the conjugacy class
$\mathsf{C}_{\lambda_{i}}$. Then $Loc_{G}(X,\{\lambda_1,\lambda_2\})$
comes equipped with a natural $(2-2d)$-shifted symplectic structure
which is a symplectic leaf of the Poisson stack $Loc_{G}(X)$.
\end{thm}

\

\noindent
Our main objective of the present paper is to understand and prove the
de Rham versions of these Betti statements.  In the next four sections
we develop the necessary framework and prove the de Rham analogue of
Theorem~\ref{thm-betti} (1). The bulk of the work goes in understanding
the algebraic $\mathcal{D}$-module theory of the formal boundary $\fb
X$ of a smooth variety where the latter is viewed as a non-commutative
space. Understanding and proving the de Rham analogue of
Theorem~\ref{thm-betti} (2) is more delicate. An important part of this
is to show the algebraicity of the derived stack of flat bundles on
$X$ which are framed by a fixed flat bundle on the formal boundary. We
prove a generic representability result for this stack of framed flat
bundles in Section~\ref{sec-represent}.  
To get a full de Rham analogue
of Theorem~\ref{thm-betti} (2) we will need to study the de Rham
versions of the gerbe twist and the strictness property for local
monodromies at intersections of components. While the first of these is
fairly straightforward, the second is quite intricate in the de Rham
context and so this discussion is left for a future work.

\section{Preliminaries}

In this section we have gathered some known and folklore results about
$\D_X$-modules on smooth varieties. We first discuss
compactness/perfection in the setting of relative $\D$-modules and its
preservation under proper push-forwards. We then recall how
$\D$-modules can be defined as graded mixed modules over the de Rham
algebra.

\subsection{Perfect relative $\D$-modules}  \label{app:relD}

\

\noindent
In this section we have gathered some basic results about $\D$-modules
in the relative setting.  Most of these results are already contained
in Gaitsgory-Rozenblyum's treatise \cite{gr}, and this part does not
claim originality. We include it here since we were unable to find a
reference treating the algebraic situation allowing for $k$ to be
non-algebraically closed, and also allowing for $\D$-modules taht are
relative over bases $Spec\, B$ with $B$ an arbitrary connective cdga.

First we discuss the compact generation and characterization of
compact objects inside quasi-coherent relative $\D$-modules. Fix $X$ a
smooth variety over $k$, and $S=Spec\, B$ an affine derived scheme. We
have $\D_X\otimes_k B$, which is a sheaf of dg-algebras over $X$. We
can therefore consider the dg-category of all sheaves of
$\D_X\otimes_k B$-modules, whose underlying $\OO_X\otimes_k B$-modules
are quasi-coherent on $X\times S$. We denote this category by
$\Dqcoh(\D_{X,B})$, and call it the dg-category of relative
$\D$-modules on $X\times S$ over $S$. An object $E\in
\Dqcoh(\D_{X,B})$ will be called \emph{\bfseries perfect} if locally
on $X$ it is given by a perfect dg-module over the dg-algebra
$\D_X\otimes_k B$.  In the special case when $B$ is a regular discrete
$k$-algebra, $\D_X\otimes_k B$ is locally a finitely generated algebra
of finite homological dimension (and thus is of finite type in the
sense of \cite[Definition~2.4]{tv}) which implies that the perfect objects are
precisely the bounded coherent $\D_X\otimes_k B$-modules. In general
the two notions do not coincide since being perfect implies in
particular being of finite tor dimension over $B$. Nevertheless we
have the following

\begin{prop}\label{p10}
The dg-category $\Dqcoh(\D_{X,B})$ is compactly generated and its
compact objects are the perfect $\D_X\otimes_k B$-modules.
\end{prop}

\noindent
\textit{Proof:} There is a forgetful functor
$$
\Dqcoh(\D_{X,B}) \longrightarrow \Dqcoh(X\times S)
$$
to the dg-category of quasi-coherent complexes on $X\times S$. This
dg-functor is conservative and continuous. Moreover, it has a left
adjoint
$$
\mathsf{ind} :
\Dqcoh(X\times S) \longrightarrow \Dqcoh(\D_{X,B})
$$ which sends a quasi-coherent complex $E$ on $X\times S$ to
$\D_X\otimes_{\OO_X}E$, with its natural $\D_X\otimes_k B$-module
structure.  It is well known that perfect complexes in $\Dqcoh(X\times
S)$ are the compact generators, and it is a formal consequence from
this that $\mathsf{ind}$-images of perfect complexes will be compact
generators of $\Dqcoh(\D_{X,B})$. These are obviously perfect
$\D_X\otimes_k B$-modules. Finally, any perfect $\D_X\otimes_k
B$-module is locally compact, and thus compact by quasi-compactness of
$X$.  \hfill $\Box$

\

\medskip

\noindent
Let now $f : X \longrightarrow Y$ be a morphism between smooth
varieties over $k$.  The usual definition gives a direct image dg-functor
$$
f_{*,B} : \Dqcoh(\D_{X,B}) \longrightarrow \Dqcoh(\D_{Y,B}).
$$ We will often drop the $B$ in the notation and simply write $f_*$.
On the level of compact generators, $f_*$ acts as follows. Let $E$ be
a perfect complex on $X\times S$, and
$\mathsf{ind}(E)=\D_{X}\otimes_{\OO_X}E$. Then we have a canonical
isomorphism $f_*(\mathsf{ind}(E))\simeq \mathsf{ind}(f_*(E))$, where
$f_*(E)$ is the direct image of $E$ as a quasi-coherent complex on
$X\times S$. In particular, when $f$ is proper the dg-functor $f_*$
preserves perfect objects. It is easy to check that the formation of
$f_*$ commutes with base change: for any morphism $B \rightarrow B'$
of connective cdga, the square
$$
\xymatrix{
  \Dqcoh(\D_{X,B}) \ar[r]^-{\otimes_B B'} \ar[d]_-{f_{*,B}} &
  \Dqcoh(\D_{X,B'}) \ar[d]^-{f_{*,B'}} \\
\Dqcoh(\D_{Y,B}) \ar[r]_-{\otimes_B B'} & \Dqcoh(\D_{Y,B'})}
$$
canonically commutes. We have thus proved the following proposition,
which is well known when $B$ is itself a smooth algebra but for which
we could not find any general reference.

\begin{prop}\label{p11}
If $f$ is proper, then $f_*$ preserves perfect objects, and its formation
commutes with change of bases $B$.
\end{prop}

\

\noindent
We recall also the following notion of holonomicity. First
recall that any coherent $\D_X \otimes_k B$-module admits a good
filtration, and that the support of the associated graded sheaf is a
well defined closed algebraic subset inside $T^{*}X \times S$.

\begin{df}\label{dhol}
Let $E \in \Dqcoh(\D_{X,B})$ be a quasi-coherent $\D_X \otimes_k
B$-module. We say that $E$ is \emph{holonomic} if it satisfies the
following two conditions.

\begin{enumerate}
\item[(1)] $E$ is perfect.
\item[(2)] There exists a conic Lagrangian algebraic subset $\Lambda
  \subset T^*X$ such that the characteristic variety of $E$ is
  contained in $\Lambda \times S$.
\end{enumerate}
\end{df}

\

\noindent
In contrast with the case of a base field, it is not true that
holonomicity for relative $\D$-modules is preserved by all six
operations.  However, this holds on a dense open subset in
$S$. Specifically, for us the following proposition will be useful.

\

\begin{prop}\label{phol}
Suppose that $B$ is a discrete noetherian $k$-algebra.  Let $E$ and
$F$ be two holonomic objects in $\Dqcoh(\D_{X,B})$. There exists a
non-empty open derived sub-scheme $Spec\, B[f^{-1}] \subset Spec\, B$
such that the tensor product $E \otimes_{\OO} F$ is a perfect
$\D_{X,B[f^{-1}]}$-module on $X \times Spec\, B[f^{-1}]$.
\end{prop}

\noindent
\textit{Proof:} Write $B_0=B_{red}$ for the reduced algebra of $B$.
Note that a given object $E \in \Dqcoh(\D_{X,B})$ is perfect if and
only if its restriction to $\Dqcoh(\D_{X,B_0})$ is perfect. Indeed, we
can use induction on the power annihilating the nil-radical of $B$ to
reduce to the case where $B$ is a square zero extension of $B_0$ by
an ideal $I$. It is easy to see that the functor sending a cdga $B$ to
the space of all quasi-coherent $\D_{X,B}$-module is 1-proximate in
the sense of formal deformation theory (see \cite{lu3}).  More
precisely, given a discrete noetherian $k$-algebra $B_0$, an ideal
$I \subset B_0$ and a derivation $d : B_0 \longrightarrow I[1]$,
consider the square-zero extension $B=B_{0}\oplus_{d} I$ classified by
$d$. The square of $\s$-categories
$$
\xymatrix{
\Dqcoh(\D_{X,B}) \ar[r] \ar[d] & \Dqcoh(\D_{X,B_0}) \ar[d] \\
\Dqcoh(\D_{X,B_0}) \ar[r] & \Dqcoh(\D_{X,B_0\oplus I[1]})}
$$
induces a full-embedding from $\Dqcoh(\D_{X,B})$ to the fiber product
of the three other terms.  This implies  that for a given $E \in
\Dqcoh(\D_{X,B})$, $E$ is a compact object if its restriction in
$\Dqcoh(\D_{X,B_0})$ is compact.

Thus it is sufficient to tackle the case where $B$ is a reduced
noetherian $k$-algebra. By picking a dense open subset in an
irreducible component we can even assume that $B$ is a smooth
domain. Let $K=Frac(B)$ be its fraction field and consider $X\times_k
K$ as a smooth variety over $K$. By the standard lore of algebraic
$\D$-modules we know that that for any smooth $K$-variety $Z$, any
smooth divisor $i_{Y} : Y \hookrightarrow Z$ given by a single
equation $f=0$ on $Z$, and any holonomic coherent $\D_{Z}$-module
$M$, there exists a Bernstein polynomial $b(M)$ for $M$ with respect to
the equation $f$. The polynomial exists as a monic polynomial over a
localization $B[f^{-1}]$ of $B$. Replacing $B$ by $B[f^{-1}]$ we can
assume that $b(M)$ exists as a monic polynomial over $B$. Now the
standard Bernstein-Kashiwara argument gives that the existence of
$b(M)$ implies that the pull-back $i_{Y}^{*}(E)$ is a bounded coherent
complex of $\D_{Y,B}$-modules with coherent cohomology, and thus is
perfect since $B$ is smooth.

Since the statement of the proposition is local on $X$, we can apply
the above reasonng to the diagonal $X \subset X \times X$, by writting
it as a complete intersection, and to the exterior tensor product
$E\boxtimes F$ of $E$ and $F$, which is manifestly a holonomic
$\D_{X\times X,B}$-module. The proposition follows.  \hfill $\Box$

\

\medskip

\subsection{Connections as graded mixed modules}

We will use freely the formalism of graded mixed $k$-modules from
\cite{ptvv,cptvv}. We denote the $\s$-category of graded mixed $k$
modules by $k-\medg$. It comes equipped with an $\s$-functor
$$
|-|:=\rh(k(0),-) : k-\medg \longrightarrow k-\dg
$$
where $k(0)$ denotes the unit in this category, i.e. the pure
weight $0$ graded mixed complex. Explicitly $|-|$ sends a graded mixed
complex $E$ to $\prod_i E(i)[-2i]$ endowed with the total differential
which is the sum of the cohomological differential and the mixed
structure.  This $\s$-functor is lax symmetric monoidal and thus
induces a corresponding $\s$-functor on algebras, modules etc.

Let $X=Spec\, A$ be a smooth affine variety over $k$ and let $\D_X$ be
the $k$-algebra of global differential operators on $X$. Consider the
de Rham algebra $\DR_X=\op{Sym}_A(\Omega_A^1[-1])$ of $X$, viewed as a
graded mixed cdga with its natural structure of a graded algebra and
with  mixed structure given by the de Rham
differential (see \cite[Section~1.1]{ptvv}). Denote by $\Dqcoh(\D_X)$ the
dg-category of complexes of left $\D_X$-modules with inverted
quasi-isomorphisms (see Section~\ref{app:relD}
for more on dg-categories of $\D$-modules).  Recall that a model for
$\Dqcoh(\D_X)$ is the dg-category of all cofibrant
$\D_X$-dg-modules. In the same way, we denote by $\DR_X-\medg$ the
dg-category of graded mixed $\DR_X$-dg-modules up to quasi-isomorphisms
(again an explicit model is the dg-category of cofibfrant graded mixed
dg-modules).  We have a natural dg-functor
$$
\DR : \Dqcoh(\D_X) \longrightarrow \DR_X-\medg,
$$
from dg-modules over $\D_X$ to graded mixed $\DR_X$-dg-modules.
The dg-functor $\DR$ is defined by
sending a (cofibrant) $\D_X$-dgmodule $E$ to its de Rham complex 
$\DR(E):=\DR_X \otimes_A E$. By definition, $\DR(E)$ is free
as a graded module over 
$\DR_X$, and its mixed structure is induced by the connection \linebreak 
$\nabla : E \longrightarrow \Omega^1_A \otimes_A E$
coming from the left $\D_X$-module structure on $E$. 

\begin{prop}\label{p1}
The dg-functor
$$\DR : \Dqcoh(\D_X) \longrightarrow \DR_X-\medg$$
is fully faithful. Its essential image consists of all objects that are
free as graded dg-modules, i.e. objects of the form $\DR_X \otimes_A E_0$ 
for some $A$-dg-module $E_0$.
\end{prop}
\textit{Proof:} To prove full faithfulness we use the following method
to compute mapping complexes inside $\DR_X-\medg$. Let $B$ be a graded
mixed cdga and $E$ and $F$ be two graded mixed $B$-dg-modules. We assume
that $E$ and $F$ are cofibrant as graded $B$-modules. Consider the
complex
$$
H(E,F):=\prod_{p\geq 0}\underline{Hom}_{B-\gdg}(E,F(p))[-p],
$$
where $F(p)$ is the graded $B$-dg-module defined by shifting the grading
by $p$ (so $\underline{Hom}_{B-\gdg}(E,F(p))$ consists of graded maps
of degree $p$).  The complex $H(E,F)$ is endowed with total
differential $D$, sending a family of elements $\{f_p\}_{p\geq 0}$ to
$$
D(\{f_p\}):=\{\nabla_F f_p + f_{p-1}\nabla_E + d(f_{p+1})\}_{p\geq 0},
$$
where $\nabla_E$ and $\nabla_F$ are the mixed structures on $E$ and
$F$, and $d$ is the cohomogical differential. Using an explicit
cofibrant model of $E$ one checks that the complex of $k$-modules
$H(E,F)$ is naturally quasi-isomorphic to the complex
$\underline{Hom}_{B-\medg}(E,F)$. This implies that the dg-functor
$\DR$ is fully faithful: for two $\D_X$-dg-modules $E$ and $F$, it
sends $\rh_{\D_X}(E,F)$ to the de Rham complex of the $\D_X$-module
$\rh_{A}(E,F)$.

For the second part of the proposition, start with a graded mixed
$\DR_X$-module $E$ which is of the form $E_0\otimes_A \DR_X$ as a
graded module. We can write $E_0$ as a filtered colimit of perfect
complexes of $A$-modules.  As the dg-functor $\DR$ is continuous and
fully faithful, it suffices to check the case where $E_0$ is
perfect. By a cell decomposition induction we can reduce to the case
where $E_0=M$ is a projective $A$-module of finite rank.  We thus have
a graded mixed $\DR_X$-module $E$ whose underlying graded module is
quasi-isomorphic to $M\otimes_A \DR_X$. We can also recover the
$\D_X$-module structure on $M$ simply by considering the map $M
\longrightarrow M\otimes_A \Omega_A^1$ induced by the mixed structure
on $E$.

This yields a canonical morphism of graded mixed dg-modules $\DR(M)
\longrightarrow E$, which by construction is a
quasi-isomorphism. \hfill $\Box$

\

\medskip

\noindent
The previous proposition extends by stackification to the case where
$X$ is a smooth scheme over $k$, or even a smooth DM-stack over
$k$. It can be stated as the existence of a full and faithful
embedding of dg-categories
$$
DR : \Dqcoh(\D_X) \hookrightarrow \DR_X-\medg,
$$
where the dg-categories $\DR_{qcoh}(\D_X)$ and $\DR_X-\medg$ are
defined by descent
$$\Dqcoh(\D_X) := \lim_{U=Spec\, A \rightarrow X}\D_U-\dg, \qquad \qquad 
\DR_X-\medg:= \lim_{U=Spec\, A \rightarrow X}\DR_U-\medg,
$$
where the limits are taken over the small \'{e}tale site of $X$ and
inside the $\s$-category of presentable dg-categories (see
\cite[Section~2]{azto}). The essential image of the dg-functor $DR$ consists of
all graded mixed $\DR_X$ dg modules which, as graded modules, are of
the form $E\otimes_{\OO_X} \DR_X$ for some quasi-coherent
$\OO_X$-module $E$.

It is also possible to extend the statement to the relative
setting. Let $B$ be a connective cdga and $X$ a smooth
DM-stack. Consider $\D_X\otimes_k B$, as a sheaf of dg-algebras, and
$\DR_X\otimes_k B$ as a sheaf of graded mixed
$B$-linear cdga (over the small \'{e}tale site of $X$).  The full
embedding $DR$ extends to a full and  faithful embedding of presentable
dg-categories
$$
DR : \Dqcoh(\D_{X,B}) \hookrightarrow (\DR_X \otimes_k B)-\medg,
$$
whose essential image consists of graded mixed modules which, as
graded modules, are of the form $E \otimes_{\OO_X} \DR_X$ for $E$ a
quasi-coherent $\OO_X \otimes_k B$-dg-module.

We conclude this part by analyzing the inverse image functor for
$\D$-modules in terms of graded mixed modules over de Rham
algebras. Let $f : X=Spec\, A' \longrightarrow Y=Spec\, A$ be a
morphism of smooth affine $k$-varieties, corresponding to a morphism
of smooth $k$-algebras $A \rightarrow A'$. We have the usual pull-back
functor of $\D$-modules
$$
f^* : \Dqcoh(\D_Y) \longrightarrow \Dqcoh(\D_X).
$$
By proposition \ref{p1} this can be seen as a dg-functor on
dg-categories of graded mixed modules which are free as graded
modules. From this point of view the functor can be described
explicitly. It is the natural functor given by base
change. Indeed, the morphism $f$ induces a morphism of graded mixed
cdga $\DR_Y \longrightarrow \DR_X$ which, in turn, defines a base
change functor on graded mixed modules. This base change is
canonically equivalent to $f^*$ when restricted to graded mixed
modules which are free as in proposition \ref{p1}.  As a final
comment, note that the above discussion also makes sense without the
affiness conditions on $X$ and $Y$, as well as in the relative setting
by tensoring with a connective cdga $B$.
 
\subsection{Graded mixed modules and equivariant objects}

We now turn to an equivalent but more conceptual description of the
dg-category of $\D$-modules, as equivariant objects inside the
dg-category of quasi-coherent modules on the shifted cotangent
stack. This will be useful later as it will allow us to reduce some
statements about $\D$-modules to statements about quasi-coherent
modules.

We let $\hH:=aut(B\mathbb{G}_a)$ be the group stack of
autoequivalences of $B\mathbb{G}_a$. It can be described explicitly as
a semi-direct product $\hH = B\mathbb{G}_a \rtimes \mathbb{G}_m$, of
$\mathbb{G}_m$ acting on $B\mathbb{G}_a$ by its natural action of
weight $1$ on $\mathbb{G}_a$. In this description, $\mathbb{G}_m$ acts
on $B\mathbb{G}_a$ by its standard action, and $B\mathbb{G}_a$ acts on
itself by translations (using the fact that $B\mathbb{G}_a$ is a
commutaive a group stack).

Recall \cite[Section~2]{azto} that there is a derived stack $\mathbb{D}g^{lp} \in
\dSt_k$ of locally presentable dg-categories with descent.  We have
the following definition.

\begin{df}\label{d1}
An \emph{\bfseries $\hH$-equivariant locally presentable dg-category} $T$
is a morphism of derived stacks $T : B\hH \longrightarrow
\mathbb{D}g^{lp}$. Locally presentable $\hH$-equivariant dg-categories
form an $\s$-category
$$
\mathbb{D}g^{lp}(B\hH):=\Map(B\hH,\mathbb{D}g^{lp}).
$$
\end{df}
Recall also that $\mathbb{D}g^{lp}$ admits a canonical extension to a
derived stack of symmetric monoidal $\s$-categories, for the tensor
product of locally presentable dg-categories of \cite{azto}. Thus we
can view a symmetric monoidal dg-category with a compatible
$\hH$-action, as a morphism $B\hH \longrightarrow
E_{\s}-Alg(\mathbb{D}g^{lp})$, from $B\hH$ to the derived stack of
$E_\s$-algebra objects in $\mathbb{D}g^{lp}$. We will not spell this
out but the interested reader can easily fill the details of this
monoidal extension.

Given an $\hH$-equivariant dg-category $T$, we can form its direct
image (see \cite{azto}) by the natural projection $p : B\hH
\longrightarrow Spec\, k$. We define the dg-category of
$\hH$-equivariant objects in $T$ to be this direct image:
$$
T^{\hH}:=p_*(T).
$$
Assume now that as in the previous section $X$ is a smooth
DM-stack, and $B$ a connective cdga. Consider $\DR_X\otimes_k B$, as a
sheaf of graded cdga on $X$, and let $(\DR_X\otimes_k B)-\dg$ be its
dg-category of (non-graded, non-mixed) dg-modules. The group $\hH$
acts on the commutative dg-algebra $(\DR_X\otimes_k B)$ in an obvious
manner: the $\mathbb{G}_m$-action is the grading and the
$B\mathbb{G}_a$-action is the mixed structure. This is formalized by
the following proposition.

\begin{prop}\label{p2}
Let $\hH$ act trivially on the dg-category $k-\dg$ of complexes of
$k$-modules. Then, there are natural equivalences of symmetric
monoidal dg-categories
$$
(k-\dg)^\hH \ \simeq  \ \Dqcoh(B\hH) \ \simeq \ k-\medg.
$$
\end{prop}
\textit{Proof:} The first equivalence holds by definition, so the
content of the proposition is the existence of the second
equivalence. For this, we let $\pi : B\hH \longrightarrow B
\mathbb{G}_m$ be the natural projection. Using this morphism we can
view $B\hH$ as an affine stack over $B\mathbb{G}_m$ whose fiber is
$K(\mathbb{G}_a,2)$.  In other words we have $B\hH \simeq
Spec_{B\mathbb{G}_m}\, A$, where $A=\pi_*(\OO_{B\hH})$ considered as an
$E_\s$-algebra in $B\mathbb{G}_m$. This algebra simply is $A=k[u]$
where $u$ is in cohomological degree $2$ and weight $1$. For any
affine stack $F=Spec\, A$, there is a symmetric monoidal $\s$-functor
$$A-Mod \longrightarrow \Dqcoh(F)$$ which makes $\Dqcoh(F)$ into the
left completion of the $A-Mod$ for the natural $t$-structure (see
\cite{lu}). This statement carries over verbatim to the relative
setting over $B\mathbb{G}_m$: there is a natural symmetric monoidal
$\s$-functor
$$A-Mod(\Dqcoh(B\mathbb{G}_m)) \longrightarrow \Dqcoh(B\hH),$$ which
is an equivalence when restricted to objects bounded on the left for
the natural $t$-structures on both sides.  Since $A=k[u]$, we have
that $\Dqcoh(B\hH)$ can be identified with the left completion of the
natural $t$-structure on the dg-category of graded
$k[u]$-dg-modules. This completion is in turn identified with the
dg-category of graded mixed complexes via the dg-functor
$$
k-\medg \longrightarrow k[u]-\dg^{gr},
$$ sending $E$ to the graded $k[u]$-module whose piece of weight $p$
is $\rh(k(p),E)$. This dg-functor is manifestly a symmetric monoidal
equivalence when restricted to graded mixed complexes which are
cohomologically bounded on the left.  This proves the
proposition. \hfill $\Box$

\

\noindent
Let $X$ be an affine smooth variety over $k$ and $B$ a connective
cdga.  $\DR_X \otimes_k B$ is a graded mixed cdga, and thus the
previous proposition can be used to view $\DR_X\otimes_k B$ as a
quasi-coherent sheaf of cdga on the stack $B\hH$.  The dg-category
$(\DR_X\otimes_k B)-\dg$ can then be seen as a natural
$E_{\s}$-algebra object in $\mathbb{D}g^{lp}(B\hH)$, or in other words
as an $\hH$-equivariant symmetric monoidal dg-category.

\begin{cor}\label{cp2}
There is a natural equivalence of symmetric monoidal dg-categories
$$
(\D_X\otimes_k B)-\dg \simeq ((\DR_X\otimes_k B)-\dg)^{\hH}.
$$ 
\end{cor}
\textit{Proof:} This follows immediately from the proposition. Indeed,
the equivalence of the proposition is symmetric monoidal, so preserves
algebras and modules over algebras.  \hfill $\Box$

\section{The formal boundary of a smooth variety} \label{sec:fb}

In this section we discuss the notion of a formal boundary of a smooth
algebraic variety $X$ over a base field $k$ of characteristic $0$. In
contrast to the Betti setting analyzed in \cite{ptbetti}, the
formal boundary does not itself exist as an algebraic variety or stack
in any form and will only be defined as a non-commutative space, i.e.
it will be defined through its category of perfect complexes.  The
requisite categories of perfect complexes have been studied recently
by several authors \cite{bete,ef,hpv}. We follow a similar approach
for the case of perfect complexes endowed with integrable connections
where many statements can be reduced to the case without
connections. However, the $\s$-category of perfect complexes with flat
connections we introduce below is a new object which can not be
recovered from the $\s$-category of perfect complexes on the formal
boundary. Thus the results of this section are new and do not follow
formally from the results of \cite{bete,ef,hpv}. 

In this section all varieties, schemes and stacks are defined over 
a base field $k$ of characteristic $0$.

\subsection{Perfect complexes on the formal boundary}

In this section we recall the notion of the formal boundary $\fb X$ of
a smooth variety $X$ studied in \cite{bete,ef,hpv}. As we prefer to
avoid any analytical aspects and constructions, we mainly follow the
approaches in \cite{ef} and \cite{hpv}.

\

\noindent
\textbf{The setting.} Let $X$ be a smooth algebraic variety. Fix an
open dense embedding $X \hookrightarrow \Xbar$ where $\Xbar$ is a
smooth and proper scheme over $k$. We moreover assume that $\Xbar$ is
chosen so that the reduced closed complement $D \subset X$ of $X$
inside $\Xbar$ is a simple normal crossing divisor on $\Xbar$.  At
some point we will also need to relax the conditions of the setting
and allow for $\Xbar$ to be a smooth and proper DM-stack, for which
the arguments are similar.  We call such an embedding $X
\hookrightarrow \Xbar$ a \emph{good compactification}.

\

For any affine scheme $Spec\, A$ with an \'{e}tale map $u : Spec\, A
\longrightarrow \Xbar$, we consider $I \subset A$ the ideal of the
pull-back $u^*(D) \subset Spec\, A$ as well as $\widehat{A}=\lim_{n}
A/I^n$ the formal completion of $A$ along $I$. When $u$ varies in the
small \'{e}tale site of $\Xbar$ we obtain a presheaf of commutative
rings on $\Xbar_{et}$, sending $u : Spec\, A \longrightarrow \Xbar$ to
$\widehat{A}$. This presheaf of commutative rings comes equipped with
a presheaf of ideals, which simply is the ideal generated by $I$
inside $\widehat{A}$.

\begin{df}\label{d2}
The \emph{\bfseries $\s$-category of perfect complexes on $\fb X$} is
defined by
$$
\iPerf(\fb X):=\lim_{Spec\, A \rightarrow \Xbar}\iPerf(Spec\, \widehat{A} - V(I)).
$$
\end{df}

\

\noindent
This definition has a version with coefficients in any derived
affine scheme $S=Spec\, B$ which goes as follows. For each $u : Spec\, A
\longrightarrow \Xbar$ in $\Xbar_{et}$ we can form the cgda
$\widehat{A\otimes B}:=\lim_{n} (A/I^n \otimes_k B)$. The ideal $I$
defines an open subset in the derived scheme $Spec\, \widehat{A\otimes
  B}$ which simply is the pull-back of $Spec\, \widehat{A}-V(I)$ by the
natural projection $Spec\, \widehat{A\otimes B} \longrightarrow Spec\,
\widehat{A}$ and we will write $Spec\, \widehat{A\otimes B} -
V(I)$ for this open derived sub-scheme. We now set
$$
\iPerf(\fb X)(S):=\lim_{Spec\, A \rightarrow \Xbar}\iPerf(Spec\,
\widehat{A\otimes B} - V(I)) \in \dgcat,
$$ and call it the {\emph\bfseries
  $\s$-category of families of perfect complexes} on
$\fb X$ parametrized by $S$. When $S$ varies in the $\s$-category of
derived affine schemes $\dAff$, $S \mapsto \iPerf(\fb X)(S)$ defines an
$\s$-functor
$$
\Perf(\fb X) : \dAff^{op} \longrightarrow \dgcat.
$$
By \cite[Proposition~3.23]{hpv} this $\s$-functor is a derived
stack for the \'{e}tale topology on $\dAff$. In the same manner, we
have the derived stack $\Perf(\widehat{\Xbar})$, of perfect complexes
on the formal completion of $\Xbar$ along $D$. For $S=Spec\, B$ i`qts
$S$-points can be defined as before
$$
\Perf(\widehat{\Xbar})(S)=\lim_{Spec\, A \rightarrow \Xbar}
\iPerf(Spec\, \widehat{A\otimes B}) \in \dgcat.
$$
Another equivalent description is as the derived mapping stack
$$\Perf(\widehat{\Xbar}) \simeq
\Map_{\dSt_k}(\widehat{\Xbar},\Perf).$$
Here the formal scheme
$\widehat{\Xbar}$ is defined as $\op{colim}_{n} \Xbar_{(n)}$, where the
colimit taken in $\dSt_k$ and \linebreak
$\Xbar_{(n)}=Spec\, \OO_X/I^{n} \subset X$ is the
$(n-1)$-th infinitesimal neighborhood of $D$ inside $\Xbar$.

\begin{df}\label{d3}
  $\Perf(\widehat{\Xbar}) \in \dSt_k$ is called the \emph{\bfseries
    derived stack of perfect complexes on $\widehat{\Xbar}$} while
  $\Perf(\fb X) \in \dSt_k$ is called the \emph{\bfseries derived
    stack of perfect complexes on $\fb X$}.
\end{df}

\

\noindent
These derived stacks admit  sheaf theoretic
interpretations.  The structure sheaf $\widehat{\OO}_{D}$ of
$\widehat{\Xbar}$ can be considered as a sheaf of commutative
$\OO_X$-algebras, sending an \'{e}tale map $Spec\, A \rightarrow \Xbar$
to the $A$-algebra $\widehat{A}$. We also have $\widehat{\OO}_{D}\simeq
\lim_n \OO_{\Xbar_{(n)}}$, where the limit is taken in the category of all
sheaves of $\OO_X$-algebras. Note that $\widehat{\OO}_D$ is in general not
a quasi-coherent sheaf on $X$. In the same manner, if $S=Spec\, B$ is
a derived affine scheme, we have a sheaf of commutative
$\OO_X$-dg-algebras $\widehat{\OO}_{D,B}$, sending an \'{e}tale map $Spec\,
A \rightarrow \Xbar$ to $\widehat{A\otimes_k B}=\lim_{n}
(A/I^n\otimes_k B)$. Again,  in general, this is not a quasi-coherent
sheaf on $X$.

Similarly, we can define a sheaf of commutative $\OO_X$-algebras
$\widehat{\OO}_D^{o}$ by locally inverting the equation of $D$ in
$\widehat{\OO}_D^{o}$. More precisely, for a derived affine scheme
$S=Spec\, B$ we send the \'{e}tale map $Spec\, A \rightarrow
\widehat{\Xbar}$ to $\Gamma(Spec\, (\widehat{A\otimes_k B}) -
V(I),\OO)$. When $Spec\, A \longrightarrow \widehat{\Xbar}$ is small
enough so that $D$ becomes principal over $Spec\, A$ (which we can
always assume for the purpose of defining the sheaf
$\widehat{\OO}_{D}^{o}$), the value of $\widehat{\OO}_{D,B}^{o}$ on
$Spec\, A \longrightarrow \widehat{\Xbar}$ is the cdga
$(\widehat{A\otimes_k B})[t^{-1}]$, where $t \in A$ is a generator of
the ideal $I \subset A$.

Both sheaves $\widehat{\OO}_{D,B}$ and $\widehat{\OO}_{D,B}^{o}$ of
cdga on $\Xbar_{et}$ are set theoretically supported on $D$, and can
therefore be considered as sheaves of cdga on the site $D_{et}$. Thus
it makes sense to consider the $\s$-categories of sheaves on $D_{et}$
which are perfect modules over the sheaves of cdgas
$\widehat{\OO}_{D,B}$ and $\widehat{\OO}_{D,B}^{o}$. Let us denote
this $\s$-categories by $\iPerf(\widehat{\OO}_{D,B})$ and
$\iPerf(\widehat{\OO}_{D,B}^{o})$. The descent result proved in
\cite[Proposition~3.23]{hpv} precisely implies that we have natural
equivalences of $\s$-categories
$$\Perf(\widehat{\Xbar})(S) \simeq \iPerf(\widehat{\OO}_{D,B}) \qquad
\Perf(\fb X)(S) \simeq \iPerf(\widehat{\OO}_{D,B}^{o})$$ which are
moreover functorial in $S=Spec\, B$

One aspect of definition \ref{d4} is that it depends a priori on a
choice of $\Xbar$.  For the perfect complexes over $\widehat{\Xbar}$
this is certainly expected but the idea is that the derived stack
$\Perf(\fb X)$ should only depend on the variety $X$. Unfortunately,
we do not know if this is the case and we could not deduce this from
the combined results of \cite{bete,ef,hpv}.  It is shown in
\cite[A.4]{hpv} (together with \cite{bete}) that the $\s$-category
$\Perf(\fb X)(k)$ of global $k$-points only depends on $X$.  However,
as noted in \cite[Appendix~A]{hpv} the setting of \cite{bete} is only for
smooth varieties and it is therefore unclear that $\Perf(\fb X)(B)$
remains independent of the choice of $\Xbar$ for a general base cdga
$B$ (already for a non-smooth commutative $k$-algebra $B$ of finite
type). To deal with this issue we introduce the full substack
$\Perf^{ex}(\fb X) \subset \Perf(\fb X)$ of \emph{\bfseries extendable
  perfect complexes} and use the categorical approach of \cite{ef}
to show that $\Perf^{ex}(\fb X)$ only depends on $X$.

To define $\Perf^{ex}(\fb X)$ consider the map of stacks in
$\s$-categories
$$\Perf(\widehat{\Xbar}) \longrightarrow \Perf(\fb X)$$
from perfect
complexes on the formal completion of $\Xbar$ along $D$ to perfect
complexes on $\fb X$. This is a morphism of stacks in stable
$\s$-categories and it therefore makes sense to define {its Karoubian
  image}. This is the substack of objects that are locally (for the
\'{e}tale topology) direct summands of objects in the essential image of
the above map. More precisely, for any affine derived scheme $S\in
\dAff$ we have a stable $\s$-functor
$$
\Perf(\widehat{\Xbar})(S) \longrightarrow \Perf(\fb X)(S),
$$
and we denote by $\Perf^{ex,pr}(\fb X)(S) \subset \Perf(\fb X)(S)$
the full sub-$\s$-category of objects that are retracts of
objects in the essential
image of $\Perf(\widehat{\Xbar})(S) \longrightarrow \Perf(\fb X)(S).$ 
When $S$ varies, this defines a full sub-prestacks $\Perf^{ex,pr}(\fb X)
\subset \Perf^{ex}(\fb X)$.

\

\begin{df}\label{d4}
The \emph{derived stack of extendable perfect complexes on $\fb X$} is
the stack associated to prestack $\Perf^{ex,pr}(\fb X)$ defined
above. It is denoted by $\Perf^{ex}(\fb X)$
\end{df}

\

Note that by definition $\Perf^{ex}(\fb X)$ is a full sub-stack in
$\Perf(\fb X)$. An important property of the stack $\Perf^{ex}(\fb X)$
is that it only depends on $X$ alone and not on the choice of $\Xbar$.

\

\begin{prop}\label{p3}
For a given $S=Spec\, B \in \dAff$, the $\s$-category $\Perf^{ex}(\fb
X)(S)$ can be reconstructed from the $k$-linear dg-category $\iPerf(X)$
of perfect complexes over the variety $X$. Moreover, this
reconstruction is functorial in $B$.
\end{prop}
\textit{Proof:} This is essentialy the main result of
\cite[Theorem~3.2]{ef} which we have bootstrapped to work over
arbitrary cdga base. First note that since $\Perf^{ex}(\fb X)$ is the
stack associated to the prestack $\Perf^{ex,pr}(\fb X)$ it is enough
to show that $\Perf^{ex}(\fb X)(S)$ can be recovered from
$\iPerf(X)$. We start with  the following lemma.

\begin{lem}\label{l1}
Let $\mathcal{K}(S)$ be the kernel of the $\s$-functor
$\Perf(\widehat{\Xbar})(S) \longrightarrow 
\Perf^{ex,pr}(\fb X)(S)$.
The sequence of stable $\s$-categories
$$\mathcal{K}(S) \hookrightarrow \Perf(\widehat{\Xbar})(S) \longrightarrow 
\Perf^{ex,pr}(\fb X)(S)$$
identifies $\Perf^{ex,pr}(\fb X)(S)$ as the
triangulated quotient of $\Perf(\widehat{\Xbar})(S)$
by $\mathcal{K}(S)$. 
\end{lem}
\textit{Proof of the lemma:} By descent, the $\s$-functor of the lemma
can be written as a finite limit over an affine cover $\mathcal{U}$ of
$\Xbar$
$$\lim_{Spec\, A \in \mathcal{U}} \iPerf(\widehat{A\otimes_k B})
\longrightarrow \lim_{Spec\, A \in \mathcal{U}}
\iPerf(\widehat{A\otimes_k B}[t^{-1}]),$$ where the affine cover
$\mathcal{U}$ has been chosen so that $D$ becomes principal on each
$Spec\, A$ and we have denoted by $t$ a local equation of $D$ in
$Spec\, A$. For a given $Spec\, A \in \mathcal{U}$ we have an exact
sequence of stable $\s$-categories
$$\mathcal{K}(S) \hookrightarrow \iPerf(\widehat{A\otimes_k B}) \longrightarrow 
\iPerf(\widehat{A\otimes_k B}[t^{-1}]).$$
But for any finite diagram of full faithful stable
$\s$-functors $T_{\alpha} \hookrightarrow T'_{\alpha}$, the induced
$\s$-functor on triangulated quotients
$$(\lim_{\alpha} T'_{\alpha}) / (\lim_{\alpha} T_{\alpha}) \longrightarrow
\lim_{\alpha} (T'_{\alpha}/T_{\alpha})$$
is fully faithful.
Therefore, the $\s$-functor $\Perf(\widehat{\Xbar})(S)/\mathcal{K}(S)
\longrightarrow \Perf^{ex,pr}(\fb X)(S)$ is always fully
faithful. Finally, by definition of extendable objects it is also
essentially surjective up to direct factors, which implies that it is
an equivalence.  \hfill $\Box$

\medskip

\

Going back to the proof of the proposition we will need a more precise
description of the kernel $\mathcal{K}(S)$. For this, we choose $K \in
\iPerf(\Xbar)$ a compact generator for $\iPerf_D(\Xbar) \subset
\iPerf(\Xbar)$, the sub dg-category of perfect complexes with supports
on $D$. The corresponding object $K \otimes_k B \in
\iPerf(\Xbar)\otimes_k B$ is a compact generator for
$\iPerf_D(\Xbar)\otimes_k B$, and this remains true after Zariski
localization on $\Xbar$: for any Zariksi open $U=Spec\, A \subset
\Xbar$, the object $K_{|U}\otimes_k B \in \iPerf(U)$ is a compact
generator for $\iPerf_D(U)$. Formal gluing for the affine $U$ (see
\cite{hpv}), tells us that we have a fibered square of dg-categories
$$
\xymatrix{ \iPerf(A\otimes B) \ar[r] \ar[d] & \iPerf(A\otimes_k
  B[t^{-1}]) \ar[d] \\ \iPerf(\widehat{A\otimes_k B}) \ar[r] &
  \iPerf(\widehat{A\otimes_k B}[t^{-1}]),}
$$
and thus an equivalence of the
kernels of the horizontal $\s$-functors. This kernel is precisely
$\iPerf_D(U)$, and thus generated by $K_{|U}\otimes_k B$. By descent, we
now have
$$ \mathcal{K}(S) \simeq \lim_{U \in \Xbar_{Zar}}\iPerf_{D\times S}(U
\times S)\simeq \iPerf_{D\times S}(\Xbar\times S).
$$
To summarize, let $C=End(K)$ be the dg-algebra of endormophisms of
the object $K$. We have an exact sequence of stable $\s$-categories
$$
\xymatrix{
  \iPerf(C\otimes_k B) \ar[r] &
  \iPerf(\widehat{\Xbar\times S}) \ar[r] & \Perf^{ex,pr}(\fb X)(S).}
$$
The $\s$-category $\iPerf(\widehat{X\times S})$ can itself be
written in terms of the dg-algebra $C$.  Again by descent we can
replace $\widehat{\Xbar}$ by an affine open sub-scheme $U=Spec\, A$
and $K$ by its restriction $K_{|U}$ to $U$. Setting $C_U :=
End(K_|U)$) we immediately see that $\iPerf(\widehat{A\otimes_k B})$
is naturally equivalent to  the dg-category $\Psi Perf(C\otimes_k B)$ of
$C_{U}$-dg-modules inside $\iPerf(B)$ (called pseudo-perfect
dg-modules relative to $B$, see \cite[Definition~2.7]{tv}). Such an equivalence is
produced by sending a perfect dg-module $E$ over $A\otimes_k B$ to
$Hom(K_{|U},E)$ as dg-module over $End(K_{|U})$. We refer to \cite{ef}
for more details.

We thus have an exact sequence of dg-categories
$$
\xymatrix{
  \iPerf(C\otimes_k B) \ar[r] & \Psi Perf(C\otimes_k B) \ar[r] &
  \Perf^{ex,pr}(\fb X)(S).}
$$
As $\iPerf(\Xbar\times S)$ is a smooth and proper dg-category over $B$,
we can now apply \cite[Theorem~3.2]{ef} to the object $K\otimes_k B
\in \iPerf(\Xbar\times S)$, which precisely states that the above
quotient can be functorially reconstructed from the $B$-linear
dg-category $\iPerf(\Xbar\times S)/\langle K\otimes_k B\rangle \simeq
\iPerf(X\times S)\simeq \iPerf(X)\otimes_k B$, and thus from $\iPerf(X)$ as
a dg-category over $k$. \hfill $\Box$

\

\begin{cor}\label{cp3}
The derived stack $\Perf^{ex}(\fb X)$ does not depend on the choice of
$\Xbar$.
\end{cor}

\

\noindent
The above corollary can be made more precise as follows. Suppose that
we have two good compactifications $\Xbar$ and $\Xbar'$ as well as a
morphism $\pi : \Xbar' \longrightarrow \Xbar$ inducing an isomorphism
over $X$. Let $\Perf(\fb X)$ and $\Perf(\fb X')$ be the two derived
stacks constructed above for $\Xbar$ and $\Xbar'$ respectively.  There
is an obvious pull-back morphism $\pi^* : \Perf(\fb X) \longrightarrow
\Perf(\fb X')$, and the corollary states that this morphism is an
equivalence of derived stacks.

Moreover, for any \'{e}tale affine $Spec\, A \longrightarrow \Xbar$, we
have a natural morphism of schemes
$$
Spec\, \widehat{A} - V(I)
\longrightarrow Spec\, A - V(I).
$$
Similary, for any $S=Spec\, B \in
\dAff$ we have a morphism of derived schemes
$$
Spec\, \widehat{A
  \otimes_k B} - V(I) \longrightarrow (Spec\, (A) - V(I))\times
S.
$$
When $A$ varies in the \'{e}tale site of $\Xbar$ and $S$ inside
derived affine schemes, we obtain by base change a natural restriction
map
$$
R : \Perf(X) \longrightarrow \Perf(\fb X),
$$
where
$\Perf(X):=\Map(X,\Perf)$ is the derived stack of perfect complexes on
$X$. Similarly, we get a restriction map
$$
R' : \Perf(X)
\longrightarrow \Perf(\fb X').
$$
Corollary \ref{cp3} and its proof then say that we have a
commutative triangle of derived stacks
$$
\xymatrix{ & \Perf(\fb X) \ar[dd]^-{\pi^*} \\
\Perf(X) \ar[ru]^-{R} \ar[rd]_-{R'} & \\
& \Perf(\fb X'), }
$$
with $\pi^*$ an equivalence. 

We do not know if the above corollary continues to hold for the bigger
stack $\Perf(\fb X)$. Because of \cite[Theorem~7.3]{hpv} the inclusion
$\Perf^{ex}(\fb X)(S) \subset \Perf(\fb X)(S)$ is an equivalence as
soon as $S$ is a smooth variety over $k$, so the restriction of
$\Perf(\fb X)$ to smooth varieties does not depend on $\Xbar$. We
believe that this remains true for a general derived affine scheme $S$
but we could not find a reference (or prove it). The question is
essentially equivalent to proving the analogue of the localization for
coherent complexes of \cite{hpv} where coherent complexes are replaced
by perfect complexes.

\subsection{Perfect complexes with flat connections on
  the formal
  boundary}

In the previous section we discussed the derived stack of perfect
complexes on the formal boundary of $X$. In this section we use
similar ideas to introduce the derived stack $\Perf^{\nabla}(\fb X)$
of perfect complexes on $\fb X$ endowed with integrable
connections. When $X=\mathbb{A}^1$ the underived version of
$\Perf^{\nabla}(\fb X)$ was extensively studied by S.~Raskin \cite{ra}
in the context of the local geometric Langlands correspondence.

We keep the setup from the previous section: we fix a smooth variety
$X$ and a good compactification $X \hookrightarrow \Xbar$, with
$D\subset \Xbar$ the divisor at infinity.  In order to define the
derived stack of perfect complexes of flat connections on $\fb X$ we
first define certain sheaves of graded mixed cgda on the small
\'{e}tale site $\Xbar_{et}$ of $\Xbar$ and then define perfect
complexes with flat connections as graded mixed dg-modules.

Let $A$ be a smooth commutative $k$-algebra of finite type.  We will
view the de Rham algebra $\DR(A)$ of $A$ over $k$, as a graded mixed
cdga over $k$. Concretely, $\DR(A)=\op{Sym}_{A}(\Omega_{A}^{1}[1])$
considered as $\mathbb{Z}$-graded cdga with zero differential and for
which $\Omega^{1}_{A}$ sits in weight $1$. The graded cdga $\DR(A)$ comes
equiped with an extra differential, namely the de Rham differential,
which we denote here by $\bepsilon$.  This additional structure makes
$\DR(A)$ into a graded mixed cdga in the sense of \cite{ptvv,cptvv}.
When $A$ is equiped with an ideal $I \subset A$, we denote by
$\widehat{\DR}(A)$ the $I$-adic completion of $\DR(A)$ which is defined by
$$
\widehat{\DR}(A):= \lim_{n} \DR(A/I^n),
$$
where the limit is taken in the category of graded mixed cdga. The
underlying graded cdga of $\widehat{\DR}(A)$ is naturally isomorphic
to $\op{Sym}_{\widehat{A}}(\widehat{\Omega}_{A}^1[1])$, the symmetric
algebra over the completion of $\Omega^1_A$. The mixed structure on
$\op{Sym}_{\widehat{A}}(\widehat{\Omega}_{A}^1[1])$ simply is the
canonical extension of the de Rham differential on $A$ to its
completion.

Let $Spec\, A \longrightarrow \Xbar$ be an \'{e}tale map, and
$I \subset A$ be the ideal of definition of the divisor $D$. We have the
completed de Rham graded mixed cdga $\widehat{\DR}(A)$. When
$Spec\, A \longrightarrow \Xbar$ varies in the small \'{e}tale site of
$\Xbar$ this defines a sheaf of graded mixed cdga $\widehat{\DR}$ on
$\Xbar_{et}$. This sheaf is set theoretically supported on $D$ and
thus defines a sheaf of graded mixed cdga on $D_{et}$. As before, this
sheaf has a version with coefficients in a cdga $B$ over $k$ denoted by
$\widehat{\DR}_{B}$.  Its values on an \'{e}tale
$U=Spec\, A \longrightarrow \Xbar$ is the graded mixed cgda
$$
\widehat{\DR}_{B}(U):= \lim_{n}(\DR(A/I^n) \otimes_k B).
$$
The sheaf $\widehat{\DR}_B$ is now a sheaf of graded mixed $B$-linear
cdga.  Note that the weight zero part of $\widehat{\DR}_B$ is the
sheaf $\widehat{\OO}_{D,B}$ constructed before. We can therefore
invert a local equation of the divisor $D$ to define
$\widehat{\DR}^{o}_B$, another sheaf of graded mixed $B$-linear
cdga. For an \'{e}tale map $U=Spec\, A \longrightarrow \Xbar$ on which
the divisor $D$ is principal with equation $t \in A$, we have
$\widehat{\DR}^{o}_B(U):=\widehat{\DR}_B(U)[t^{-1}]$.  The part of
weight zero in $\widehat{\DR}^{o}_B(U)$ is of course the sheaf
$\widehat{\OO}_{D,B}^{o}$ defined in the previous section.

For $S=Spec\, B \in \dAff_k$, we let $\iPerf^{\nabla}(\fb X)(S)$ be the
dg-category of sheaves $E$ of graded mixed
$\widehat{\DR}^{o}_B(U)$-dg-modules which are locally free of weight $0$
in the following sense: locally on $\Xbar_{et}$, the underlying graded
$\widehat{\DR}^{o}_B$-dg-module $E$ (obtained by forgetting the mixed
structure) is of the form
$\widehat{\DR}^{o}_B\otimes_{\widehat{\OO}_{D,B}^{o}}E(0)$ for some perfect
$\widehat{\OO}_{D,B}^{o}$-module $E(0)$ of weight zero.
When $S=Spec\, B$ varies in $\dAff$, the dg-categories
$\iPerf(\fb X)(S)$ define an dg-functor
$\Perf^{\nabla}(\fb X) : \dAff^{op} \longrightarrow \dgcat$.  There is
an obvious forgetful map of derived prestacks
$$\Perf^{\nabla}(\fb X) \longrightarrow \Perf(\fb X)$$
sending a graded mixed dg-module to its part of weight $0$.

\

\begin{df}\label{d5}
  $\Perf^{\nabla}(\fb X)$ is called the \emph{\bfseries derived
    pre-stack of perfect complexes with flat connections on $\fb X$}.
  The \emph{\bfseries derived pre-stack of extendable perfect
    complexes with flat connections on $\fb X$} is defined to be the
  fiber product of derived pre-stacks
  $$
  \Perf^{\nabla,ex}(\fb X) \times_{\Perf(\fb X)}\Perf^{ex}(\fb X).
  $$
\end{df}

\

\noindent
By construction, $\Perf^{\nabla,ex}(\fb X)$ is a full derived
sub-prestack in $\Perf^{\nabla}(\fb X)$ defined by the local condition
"the underlying perfect complex is extendable". The main result of
this section is the following descent and invariance statement.

\

\begin{prop}\label{p4}
With the notations above we have:
\begin{enumerate}
\item The derived pre-stacks $\Perf^{\nabla}(\fb X)$ and
  $\Perf^{\nabla,ex}(\fb X)$ are stacks.
\item The derived stack $\Perf^{\nabla,ex}(\fb X)$ only depends on
  $X$.
\end{enumerate}
\end{prop}
\textit{Proof:} The key to the proof of this proposition is the
interpretation of graded mixed structures as actions of the group
stack $\hH:=B\mathbb{G}_a \rtimes \mathbb{G}_m$ (see
Proposition~\ref{p1}).  For a graded mixed cdga $\Omega$, the group
stack $\mathcal{H}$ acts on $\Omega$ by cdga automorphisms, where the
$\mathbb{G}_m$-action provides the grading and the
$B\mathbb{G}_a$-action induces the mixed structure. This action
induces an action of $\mathcal{H}$ on the $k$-linear dg-category
$\iPerf(\Omega)$ of perfect dg-modules over $\Omega$. The dg-category of
graded mixed $\Omega$-modules which are perfect as $\Omega$-dg-modules
can be recovered by taking invariants (see \ref{p1})
$$
\iPerf^{gr,\bepsilon}(\Omega) \simeq \iPerf(\Omega)^{\mathcal{H}}.
$$
This presentation of graded mixed dg-modules implies the statement of
the proposition as follows.

For $(1)$, the derived prestack $\Perf^{\nabla}(\fb X)$ is obtained as
follows. We start with the prestack $\iPerf(\widehat{\DR}^{o})$ of
perfect $\widehat{\DR}^{o}$-dg-modules, where $\widehat{\DR}^{o}$ is
simply considered as a sheaf of graded cdga. This is a derived
prestack with values in $\mathcal{H}$-equivariant dg-categories. It is
moreover a stack, which follows by noticing that $\widehat{\DR}^{o}$ is
a cdga inside $\Perf(\fb X)$ and by using
\cite[Proposition~3.23]{hpv}.  This implies that its fixed points by
$\mathcal{H}$ remain a stack (because taking fixed points commutes
with taking limits). This stack is denoted by
$\iPerf^{gr,\bepsilon}(\widehat{\DR}^{o})$, and is the stack of graded
mixed $\widehat{\DR}^{o}$-dg-modules which are perfect over
$\widehat{\DR}^{o}$. But $\Perf^{\nabla}(\fb X)$ is a sub-prestack of
the stack $\iPerf^{gr,\bepsilon}(\widehat{\DR}^{o})$ which is defined by a
local condition and thus is a stack. The fact that
$\Perf^{\nabla,ex}(\fb X)$ is also a stack now follows from the fact
that it is defined as a fiber product of stacks.

For $(2)$ we use a similar argument. The derived stack
$\Perf^{\nabla,ex}(\fb X)$ can be expressed as a full sub-stack of the
fixed points by $\mathcal{H}$ acting on $\widehat{\DR}^{o}$-dg-modules
inside $\Perf^{ex}(\fb X)$ (note that as a graded cdga
$\widehat{\DR}^{o}$ lives in $\Perf^{ex}(\fb X)$).  Therefore, to prove
that $\Perf^{\nabla,ex}(\fb X)$ is independant of the choice of
$\Xbar$ we have to check that the stack of $\mathcal{H}$-equivariant
dg-categories $\iPerf(\widehat{\DR}^{o})$ only depends on $X$. This
reduces to the following lemma.

\begin{lem}\label{l2}
Let $\pi : \Xbar' \longrightarrow \Xbar$ be a morphism between two
good compactifications of $X$.  Let $D'=\pi^{-1}(D)$ so that $\pi$
induces an isomorphism between $\Xbar'-D'$ and $\Xbar-D$.  Let
$\widehat{\DR}_{\Xbar}^{o}$ and $\widehat{\DR}_{\Xbar'}^{o}$ be the
corresponding  sheaves of graded mixed cdga constructed above.
Then, for any $Spec\, B \in \dAff$, we have.

\begin{enumerate}

\item There is a pull-back map 
  $ f_{\pi} : \pi^{-1}(\widehat{\DR}_{\Xbar,B}^{o}) \longrightarrow
  \widehat{\DR}_{\Xbar',B}^{o}
  $ of sheaves of graded mixed cdga on $\Xbar'_{et}$.

\item The map $f_{\pi}$
  induces an equivalence of dg-categories
$\pi^* : \iPerf(\widehat{\DR}_{\Xbar,B}^{o}) \simeq
  \iPerf(\widehat{\DR}_{\Xbar',B}^{o}).$
\end{enumerate}
\end{lem}

Before giving a proof of the lemma, let us explain how this finishes
the proof of the proposition.  The fact that $f_\pi$ exists implies
that the dg-functor $\pi^*$ also exists by simply pulling back graded
mixed dg-modules. Moreover, as $f_\pi$ is a morphism of graded mixed
cdga, it is clear that the dg-functor $\pi^*$ is naturally
$\mathcal{H}$-equivariant. As it is an equivalence it also induces an
equivalence on the fixed points dg-categories, and the result follows
immediately by considering the full sub-dg-categories corresponding to
$\Perf^{\nabla,ex}(\fb X)$.

\

\noindent
\textit{Proof of the lemma:} The existence of the map $f_\pi$ simply
follows from the fact that the assignments  $A \mapsto \DR_B(A)$, \
$A \mapsto \widehat{\DR}_B(A)$, \  and $ A \mapsto
\widehat{\DR}_B(A)[t^{-1}]$ are  functors from smooth $k$-algebras of
finite type to graded mixed cdga.  To prove $(2)$, we observe that
$\widehat{\DR}_{\Xbar,B}^{o}$ and $\widehat{\DR}_{\Xbar',B}^{o}$, when
considered as sheaves of cdga are perfect over $\widehat{\OO}_{D,B}^{o}$
and $\widehat{\OO}_{D',B}^{o}$ and extendable. They can thus be
considered as graded cdga inside the symmetric monoidal the
dg-categories $\Perf^{ex}(\fb X)(B)$ and $\Perf^{ex}(\fb X')(B)$. By
corollary \ref{cp3} we know that pull-back along $\pi$ induces an
equivalence of symmetric monoidal dg-categories
$$\pi^* : \Perf^{ex}(\fb X)(B) \simeq \Perf^{ex}(\fb X')(B).$$ To
finish the proof it remains to show that the symmetric monoidal
equivalence $\pi^*$ sends the cdga $\widehat{\DR}_{\Xbar,B}^{o}$ to
$\widehat{\DR}_{\Xbar',B}^{o}$.  There are canonical restriction maps
$$R : \Perf(X) \longrightarrow \Perf(\fb X) \qquad R' : \Perf(X)
\longrightarrow \Perf(\fb X')$$ and we have $\pi^*\circ R' \simeq
R$. Moreover, by construction we have
$$\widehat{\DR}_{\Xbar,B}^{o} \simeq R(\DR_X) \qquad
\widehat{\DR}_{\Xbar',B}^{o} \simeq R'(\DR_X),$$ where
$\DR_X=\op{Sym}_{\OO_X}(\Omega^1_X[1])$ as a sheaf of perfect cdga
over $\OO_X$.  This completes the proof of the lemma and of
proposition \ref{p4}. \hfill $\Box$

\

\noindent
To finish this section, note that as for the case of perfect
complexes, there is a restriction morphism
$$
\mathsf{R} : \Perf^{\nabla}(X) \longrightarrow \Perf^{\nabla,ex}(\fb X)
\subset \Perf^{\nabla}(\fb X),
$$
from the derived stack $\Perf^{\nabla}(X)$ of perfect
complexes on $X$ endowed with flat connections to the derived stack of
extendable perfect complexes with flat connections on the formal boundary of
$X$. It is defined as follows. First of all the derived stack
$\Perf^{\nabla}(X)$ is defined as the derived stack of graded mixed
dg-modules over $\DR_X$, the de Rham algebra of $X$, which are perfect
of weight $0$. More precisely, for $S=Spec\, B \in \dAff$, then
$\Perf^{\nabla}(X)(S)$ is defined to be the $\s$-category of graded
mixed $\DR_X \otimes_k B$-dg-modules $E$, such that $E\simeq E(0)
\otimes_{\OO_X}\DR_X$ as a graded dg-modules over $B$, and where
$E(0)$ is perfect over $\OO_X \otimes_k B$. The restriction map
$\mathsf{R}$ is then induced by the natural morphism of sheaves of
graded mixed cdga over $\Xbar_{et}$
$$\DR_X \otimes_k B \longrightarrow \widehat{\DR}_{B}^{o}.$$
Locally on an \'{e}tale affine $Spec\, A \longrightarrow \Xbar$ on which 
$D$ is principal with equation $t \in A$, this morphism is the natural map
$$
\DR(A \otimes_k B)[t^{-1}] \longrightarrow
\widehat{\DR}(A\otimes_k B)[t^{-1}]
$$
induced by the completion morphism $A \otimes_k B \longrightarrow
\lim_n (A/I^n \otimes_k B)$. 

This defines a restriction map $\mathsf{R} : \Perf^{\nabla}(X)
\longrightarrow \Perf^{\nabla}(\fb X)$, which covers the restriction
map of perfect complexes $R : \Perf(X) \longrightarrow \Perf(\fb
X)$. As the later map factors through extendable perfect complexes
(because any perfect complex on $X \times S$ extends to $\Xbar\times
S$ up to a retract), we also get a restriction map $\Perf^{\nabla}(X)
\longrightarrow \Perf^{\nabla,ex}(\fb X)$,

\begin{df}\label{d6}
The \emph{\bfseries de Rham restriction map} is the morphism of
derived stacks
$$
\mathsf{R} : \Perf^{\nabla}(X) \longrightarrow \Perf^{\nabla,ex}(\fb X)
$$
defined above. 
\end{df}

\subsection{De Rham cohomology of the formal boundary 
and compactly
  supported cohomology} \label{s-cptsup}

To finish this section, let us describe the $Hom$-complexes of the
dg-category $\Perf^{\nabla}(\fb X)(B)$ in terms of hypercohomology of
complexes of sheaves on $D$ and relate this to the notion of compactly
supported de Rham cohomology. The notion of de Rham cohomology with
compact supports already appeared in \cite{bacafi}, but our treatment
here is new as it is based on the theory of Tate objects and their
duality (see \cite{he}), which makes the theory also available over
any base cdga $B$. In this part we give the constructions and
definitions of compactly supported cohomology. The duality itself
is studied  in section~\ref{s-orientation}.

Again we fix a good compactification $X\hookrightarrow \Xbar$ with
divisor at infinity $D$. For any connective cdga $B$ and any object
$E \in \Perf^{\nabla}(\fb X)(B)$, we define a sheaf of
$B$-dg-modules on $D_{Zar}$ as follows. By definition, $E$ is a sheaf
of graded mixed modules over $\widehat{\DR}_{B}^{o}$. We define $|E|$ to
be the sheaf of $B$-modules $\underline{Hom}_{\medg}(k,E)$, of graded
mixed morphisms from the unit $k$ to $E$. Note that a priori $|E|$ is 
given as an infinite product
$$
|E|=\prod_{i\geq 0}E(i)[-2i],
$$
where the differential is the sum of the cohomological differential
and the mixed structure. However, in our situation this infinite
product is in fact a finite product, as $E(i)$ is non-zero only for a
finite number of indices $i$ (because $E$ is free as a graded module
and $\widehat{\DR}_{B}^{o}$ only has weights in the interval$[0,d]$
where $d=dim_{k} X$).  We will call the sheaf of $B$-dg-modules $|E|$
the \emph{\bfseries de Rham complex of $E$ completed along $D$}.  With
this notation we now have the following

\begin{df}\label{d7}
The \emph{\bfseries  de Rham cohomology of $\fb X$ with coefficients
in $E$} is the $B$-module defined by
$$\HH_{DR}(\fb X,E):=\HH(D,|E|) \in B-\dg.$$
\end{df}

Going back to the problem of computing $Hom$-complexes let $E$ and $F$
be two objects in $\Perf^{\nabla}(\fb X)(B)$.  The dg-category of
graded mixed modules over $\widehat{\DR}_{B}^{o}$ has a canonical
symmetric monoidal structure, for which the tensor product is given by
tensoring the underlying $B$-dg-modules (see \cite{ptvv,cptvv}). Since
perfect complexes of $\widehat{\OO}_{D,B}^{o}$-modules form a rigid
symmetric monoidal dg-category it follows that
$\Perf^{\nabla}(\fb X)(B)$ is also rigid. We can then form
$E^{\vee}\otimes_{\widehat{\OO}_{D,B}^{o}} F$, which is a new graded
mixed $\widehat{\DR}_{B}^{o}$-module and an object in
$\Perf^{\nabla}(\fb X)(B)$. To simplify notation we will
simply write $E^{\vee}\otimes F$ for this object.  We then have a natural
quasi-isomorphism
$$
\underline{Hom}_{\Perf^{\nabla}(\fb X)(B)}(E,F)\simeq \HH_{DR}(\fb X,E^\vee
\otimes F),
$$
giving us the desired interpretation of mapping
complexes of $\Perf^{\nabla}(\fb X)(B)$ in terms of de Rham cohomology
of $\fb X$.

\begin{rmk}
  For any connective $B$, and any object
  $E \in \Perf^{\nabla}(\fb X)(B)$ the complex of sheaves $|E|$ on
  $D_{Zar}$ is built out of acyclic sheaves on affines. Therefore the
  hyper-cohomology complex $\mathbb{H}(D,|E|)$ can be computed by a
  finite limit using an affine cover of $D$.  In particular, if $|E|$
  is locally perfect as a $B$-module on $D$, then $\mathbb{H}(D,|E|)$
  is a perfect $B$-module.
\end{rmk}

\

\noindent
We now show how the formal boundary $\fb X$ can be used to define a
notion of cohomology with compact supports, both for perfect complexes
and perfect complexes with flat connections on $X$.  We start with a
connective cdga $B$ and a perfect complex $E$ on $X\times S$, where
$S=Spec\, B$. As explained in the previous section (right before definition
\ref{d6}) we have its restriction $R(E) \in \Perf(\fb X)(B)$, and by
functoriality an induced map on cohomology
\begin{equation} \label{eq-R.on.cohomology}
  \HH(X,E)=\underline{Hom}(\OO_X,E) \longrightarrow \HH(\fb X,R(E))=
  \underline{Hom}(R(\OO_X),R(E)).
\end{equation}

The \emph{\bfseries cohomology of $X$ with compact supports and with
  coefficients in $E$} is defined to be the homotopy fiber of the map of
complexes \eqref{eq-R.on.cohomology}. It
is denoted by
$$\HH_{c}(X,E):=\mathsf{fib}\left(\HH(X,E) \longrightarrow
  \HH(\fb X,R(E))\right) \in B-\dg.
$$
By construction this is a $B$-dg-module. This is not quite enough for
our purpose, as this $B$-dg-module turns out to be the realization of
a natural pro-object that we will now describe. This pro-structure is
going to be very important for us, as it will allow us to work with
compactly supported cohomology as dual of cohomology, even if the
later is infinite dimensional. For simplicity we assume
that $E$ extends to our fixed good compactification as a perfect
complex $\barE$ on $\Xbar\times S$.  It is not always possible to find
$\barE$ in general although it always exists if $B=k$ (because
$K_{-1}(X)=0$). Moreover, such an extension always exists up to a
retract, so assuming the existence of $\barE$ is thus not a real
restriction.  Note also that in all our applications $E$ will always
come with an extension to $\Xbar$.

Using the formal gluing formalism of \cite{hpv}, we obtain 
a cartesian square of $B$-dg-modules
$$
\xymatrix{
\HH(\Xbar,\barE) \ar[r] \ar[d] & \HH(X,E) \ar[d] \\
\HH(\widehat{\Xbar},\widehat{\barE}) \ar[r] & \HH(\fb X,R(E)).
}
$$
Here $\widehat{\Xbar}$ is the formal completion of $X$ along $D$, and 
$\HH(\widehat{\Xbar},\widehat{\barE})$ is defined by
$$
\HH(\widehat{\Xbar},\widehat{\barE})
:=
\lim_{n} \HH\left(\Xbar_{(n)}, \, \jfrak_{n}^{*}\barE\right),
$$
where $\jfrak_n : \Xbar_{(n)}:= Spec\ (\OO_{\Xbar}/I_D^n)
\longrightarrow \Xbar$
is the $(n-1)$-th infinitesimal thickening of $D$ inside $\Xbar$.

From the diagram above we 
have that $\HH_c(X,E)$ can also be described as the fiber of  the map
\begin{equation} \label{eq:HHc}
\HH(\Xbar,\barE) \longrightarrow \lim_{n}
\HH(\Xbar_{(n)},\jfrak^{*}_{n}\barE).
\end{equation}
The morphism \eqref{eq:HHc} can itself be considered as a morphism of
pro-objects in $\iPerf(B)$. This allows us to define a
pro-perfect $B$-module by
$$
\HHc(X,E):= \mathsf{fib}\left(
  \HH(\Xbar,\barE) \longrightarrow \text{``$\lim_{n}$''}
  \HH(\Xbar_{(n)},\jfrak^{*}_{n}\barE)\right)
\in \mathit{ProPerf}(B).
$$
It is easy to show that this definition does not depend on choosing 
either $\Xbar$ or $\barE$, but we will not do it here. For us this
will be a consequence of Serre duality with supports which is studied
in section~\ref{s-orientation}, as the dual $B$-module turns out to be
canonically equivalent to $\HH(X,E^{\vee}\otimes_{\OO_X} \omega_X)$,
which only depends on $X$ and $E$. For future reference we record the
following 

\begin{df}\label{d7'}
 The \emph{\bfseries refined cohomology with compact
    supports} of $X$ with coefficients in $E$, is the pro-perfect
  $B$-module $\HHc(X,E)$ defined above.
\end{df}

One nice aspect of the refined version of compactly supported
cohomology is that it is manifestly compatible with base changes of
$B$.  Let $B \rightarrow B'$ be any morphism of connective cdga,
then the natural map
$$\HHc(X,E) \widehat{\otimes}_B B' \longrightarrow \HHc(X,E\otimes_B B')$$
is an equivalence of pro-perfect $B'$-modules. Here we have denoted by 
$$
\widehat{\otimes}_B B' : \mathit{ProPerf}(B) \longrightarrow
\mathit{ProPerf}(B')
$$
the functor induced on pro-objects by the usual base change
$\otimes_B B' : \iPerf(B) \longrightarrow \iPerf(B')$.

Another feature  of the  refined compactly supported
cohomology  comes  from the existence of a fiber sequence of $B$-modules
$$
\xymatrix{\HH_c(X,E)  \ar[r] & \HH(X,E) \ar[r] & \HH(\fb X,R(E)).}
$$
The first map in this sequence arises from  a natural morphism of
ind-pro-perfect $B$-modules \linebreak
$\HHc(X,E) \rightarrow \HH(X,E)$, where
$\HH(X,E)$ is considered as an ind-perfect $B$-module via the canonical
equivalence  $\mathit{IndPerf}(B)\simeq B-\dg$.  This implies that
$\HH(\fb X,R(E))$ is itself the realization of an ind-pro-perfect
module $\widetilde{\HH}(\fb X,R(E)) \in \mathit{IndProPerf}(B)$,
sitting in a triangle
$$
\xymatrix{\HHc(X,E)  \ar[r] & \HH(X,E) \ar[r] & \widetilde{\HH}(\fb
  X,R(E)).}
$$
By construction the ind-pro-perfect object $\widetilde{\HH}(\fb X,R(E))$
is an extension of 
a pro-perfect by an ind-perfect, and thus by definition is a Tate
$B$-module in the sense of \cite{he}. 

We now turn to the case of an object $E \in \Perf^{\nabla}(X)(B)$. The
naive de Rham cohomology of $E$ with compact supports is again defined as
$$
\HH_{c,DR}(X,E):=
\mathsf{fib}\left(\HH_{DR}(X,E) \longrightarrow \HH_{DR}(\fb
X,\mathsf{R}(E))\right) \in B-\dg.
$$
As before this $B$-module is the realization of a natural
pro-perfect $B$-module denoted by $\HHcd(X,E)$.  We assume again that
the underlying perfect complex $E(0)$ of $E$ extends to a perfect
complex $\barE(0)$ on $\Xbar\times S$.  Using this, one immediately
checks that the sheaf of $B$-modules $|R(E)|$ on $D$ has a natural
structure of a sheaf of ind-pro $B$-modules. Indeed, it is of the form
$\oplus_i R(E(i))[-2i]$ with a suitable differential.  Each $R(E(i))$
is itself of the form $\barE(0) \otimes_{\OO_{\Xbar}}\Omega^i_{\Xbar}
\otimes_{\OO_{\Xbar}}\widehat{\OO}_{D,B}^{o}$.  As
$\widehat{\OO}_{X,B}^0$ has a canonical ind-pro structure, and since
the functor $\barE(0) \otimes_{\OO_X} \Omega_X^i \otimes_{\OO_X}(-)$
commutes with limits and colimits of $\OO_X$-modules, we see that each
$R(E)(i)$ is the realization of a canonical sheaf of ind-pro
$B$-modules. Moreover, since $\hat{\OO}_{D,B}^{o}$ is ind-pro perfect
as a $B$-module, this endows $|E|$ with a natural structure of sheaf
of Tate $B$-modules.  This provides a canonical Tate structure on the
hyper-cohomology of $D$ with coefficients in $|E|$, that is of the de
Rham cohomology of $\fb X$ with coefficients in $\mathsf{R}(E)$. We
denote this Tate $B$-module by $\widetilde{\HH}_{DR}(\fb
X,\mathsf{R}(E))$.

The restriction map $\mathsf{R}$ induces a morphism
$\HH_{DR}(X,E) \longrightarrow \widetilde{\HH}_{DR}(\fb X,\mathsf{R}(E))$,
which is a morphism of ind-pro perfect $B$-modules if one endows the
left hand side with the canonical structure of an ind-perfect
$B$-module. It thus lifts to a morphism of Tate $B$-modules
$$
\widetilde{\HH}_{DR}(X,E) \longrightarrow \widetilde{\HH}_{DR}(\fb
X,\mathsf{R}(E)).
$$
With this notation we can now formulate the following 

\begin{df}\label{d8}
 The \emph{\bfseries refined de Rham
    cohomology of $X$ with compact supports with coefficients in $E$}
  is the Tate $B$-module defined by
  $$\HHcd(X,E):= \mathsf{fib}\left(\widetilde{\HH}_{DR}(X,E)
  \longrightarrow \widetilde{\HH}_{DR}(\fb X,\mathsf{R}(E))\right).$$
\end{df}

\

\noindent
It is instructive to note that the ind-pro-perfect $B$-module
structure on $\HHcd(X,E)$ is in fact pro-perfect (in particular it is
a Tate $B$-module in the sense of \cite{he}). This can be seen by
reducing to the previously treated case of perfect complexes without
connections. Indeed, the complexes of sheaves $|E|$ and
$|\mathsf{R}(E)|$ are canonically filtered using their Hodge
filtrations.  The graded pieces of the Hodge filtration on
$\HHcd(X,E)$ are $\HHc(X,\Omega_{X}^i\otimes_{\OO_X}E(0))[-i]$, and
thus are pro-perfect.  Since this filtration is finite, we deduce that
ind-pro object $\HHcd(X,E)$ is filtered with associated graded being
pro-perfect. This implies that $\HHcd(X,E)$ itself is pro-perfect.  We
thus have proven the following corollary.

\begin{cor}\label{c4}
  The ind-pro perfect $B$-module $\HHcd(X,E)$ is
  pro-perfect. Furthermore the ind-pro perfect $B$-module
  $\widetilde{\HH}_{DR}(\fb X,\mathsf{R}(E))$ is a
  Tate $B$-module in the sense
  of \cite{he}.
\end{cor}

\

\noindent
As in the case of perfect complexes, the formation of $\HHcd(X,E)$
commutes with base change over $B$: for any $B \rightarrow B'$ of
connective cdga, the natural morphism
$$\HHcd(X,E)\widehat{\otimes}_B B' \longrightarrow \HHcd(X,E\otimes_B B')$$
is an equivalence of pro-perfect $B'$-modules. 

\section{Formal properties of moduli functors}

We start by recalling some of the general formal properties of derived
stacks (see \cite{hagII}, \cite{lu2}).  Let $F \in \dAff_k$ be a
derived stack over $k$. For any derived affine scheme
$u : U =Spec\, B \longrightarrow F$ mapping to $F$, and any connective
$B$-dg-module $M$, we can define the space of derivations of $F$ on
$U$ with coefficients in $M$, as the fiber at $u$ of the
restriction map
$$
F(B\oplus M )\longrightarrow F(B),
$$
where $B\oplus M$ is the trivial square zero extension of $B$ by
$M$. Denote this space by $Der_u(F,M) \in \T$. For any morphism
$B \longrightarrow B'$ of connective cdga and any connective
$B'$-dg-module $M'$, we have a canonical morphism
$B\oplus M' \longrightarrow B' \oplus M'$ covering the map
$B \rightarrow B'$. Therefore, for any commutative diagram of derived
stacks
$$\xymatrix{
U=Spec\, \ar[dr]_-{u}\ar[rr]^-{f} & & U'=Spec\ ,B' \ar[dl]^-{u'} \\
 & F & }$$
 there is a natural induced morphism on the corresponding spaces of derivations
$$f^* : Der_u(F,M) \longrightarrow Der_{u'}(F,M').$$

\begin{df}\label{d9} \
\begin{enumerate}
\item[(1)] The derived stack $F$ has a \emph{\bfseries cotangent complex at
    $u : U=Spec\, B \longrightarrow F$} if there is an eventually
  connective $B$-dg-module $\mathbb{L}_{F,u}$ and functorial
  equivalences
$$Map_{B-mod}(\mathbb{L}_{u,F},M) \simeq Der_u(F,M).$$
\item[(2)] We say that $F$ has a (gobal) \emph{\bfseries cotangent complex}
  if it has cotangent complexes at all maps
  $u : U=Spec\, B \longrightarrow F$, and if moreover for any commutative
  diagram
$$\xymatrix{
U=Spec\, \ar[dr]_-{u}\ar[rr]^-{f} & & U'=Spec\ ,B' \ar[dl]^-{u'} \\
 & F & }$$
the induced morphism $Der_u(F,M) \rightarrow Der_{u'}(F,M')$
is an equivalence.
\end{enumerate}
\end{df}

\

\noindent
It is shown in \cite{hagII,lu2} that $\mathbb{L}_{F,u}$, if it exists,
is uniquely characterized by the $\s$-functor $Der_{u}(F,-)$. Also,
condition $(2)$ can be reformulated as the statement that the natural
morphism
$\mathbb{L}_{u,F} \otimes_{B}B' \rightarrow \mathbb{L}_{u',F}$ is an
equivalence of dg-modules.

Let $Spec\, B \in \dAff_k$ be a derived affine scheme and $M$ a
connective $B$-module.  Let $d : B \longrightarrow M[1]$ be a
$k$-linear derivation, which by definition means a section of
$B \oplus M \longrightarrow B$ inside the $\s$-category of cdga over
$k$.  Recall that the square zero extension of $B$ by $M$ with
respect to $d$, denoted by $B\oplus_d M $ is defined by the cartesian
square of cdga (see \cite{hagII})
$$\xymatrix{
B\oplus_d M \ar[r] \ar[d] & B \ar[d]^-{0} \\
B\ar[r]_-{d} & B\oplus M[1]
}$$
where $0$ denotes the natural inclusion of $B$ as a direct factor in the trivial
square zero extension $B\oplus M[1]$. 

\begin{df}\label{d10}
Let $F$ be a derived stack. 
\begin{enumerate}

\item We say that $F$ is \emph{\bfseries
    inf-cartesian} if for any $B$, $M$ and $d$ as above
the square
$$\xymatrix{F(B\oplus_d M) \ar[r] \ar[d] & F(B) \ar[d]^-{0} \\
F(B)\ar[r]_-{d} & F(B\oplus M[1])}$$
is cartesian.

\item We say that $F$ is \emph{\bfseries nil-complete} if for any
  $Spec\, B \in \dAff_k$
with Postnikov tower $\{B_{\leq n}\}_{n}$ the natural morphism
$$F(B) \longrightarrow \lim_n F(B_{\leq n})$$
is an equivalence.
\end{enumerate}
\end{df}

\

\noindent
Suppose now that $F$ is a derived stack which is inf-cartesian. For
any $x : Spec\, B \longrightarrow F$, we have an $\s$-functor
$$\mathbb{T}_{F,x} : B-Mod^{c} \longrightarrow \T,$$
from connective $B$-modules to spaces, that sends $M$ to the fiber of
$F(B\oplus M) \longrightarrow F(B)$ at the point $x$. This
$\s$-functor restricts to the full sub-$\s$-category of $B$-modules of
the form $B[i]^n$ for various $i\geq 0$ and various $n$. Because $F$
is inf-cartesian, the $\s$-functor $\mathbb{T}_{F,x}$ preserves finite
products as well as the looping construction $\Omega_*$ (i.e. the
natural map
$\mathbb{T}_{F,x}(M[-1]) \rightarrow \Omega_*(\mathbb{T}_{F,x}(M))$ is
an equivalence of spaces).  This implies that there exists a unique
$B$-dg-module $\mathbb{T}_{F,x}$ such that
$\mathbb{T}_{F,x}(B[i]^n) \simeq
Map_{B-Mod}(B[-i]^n,\mathbb{T}_{F,x})$ for all $i\geq 0$ and $n$. We
still denote this complex by $\mathbb{T}_{F,x}$, and call it the
\emph{\bfseries tangent complex of $F$ at $x$}.  The following result
is an easy criterion for existence of cotangent complexes.

\begin{lem}\label{l3}
  Let $F$ be a derived stack which is inf-cartesian and
  $x : Spec\, B \longrightarrow F$. 
Assume that the two conditions below are satisfied.
\begin{enumerate}
\item[(1)] The $\s$-functor $M \mapsto \mathbb{T}_{F,x}(M)$ commutes with
  arbitrary colimits.
\item[(2)] The $B$-module $\mathbb{T}_{F,x}$ is perfect.
\end{enumerate}
Then 
$F$ has a cotangent complex $\mathbb{L}_{F,x}$ at $x$ and moreover we have
$\mathbb{L}_{F,x}$ is naturally identified with $\mathbb{T}_{F,x}^{\vee}$, 
the $B$-linear dual of $\mathbb{T}_{F,x}$.
\end{lem}

\noindent
\textit{Proof:} We consider the two $\s$-functors
$$B-Mod^c \longrightarrow \T,$$
sending $M$ to either $Map_{B-Mod}(B,\mathbb{T}_{F,x}\otimes_B M)$ or 
$\mathbb{T}_{F,x}(M)$. There is a canonical equivalence of functors
\begin{equation} \label{eq-tangent}
  Map_{B-Mod}(B,\mathbb{T}_{F,x}\otimes_B -) \ \cong \  \mathbb{T}_{F,x}(-)
\end{equation}
when these functors are restricted to the full
sub-$\s$-category of objects of the form $B[i]^n$.  However, these
objects generate $B-Mod^c$ by colimits, so by condition $(1)$ the map
\eqref{eq-tangent}  extends to  an equivalence of $\s$-functors defined on the
whole $\s$-category $B-Mod^c$. In formulas - for any connective
$B$-module $M$ we have a natural equivalence
$$Map_{B-Mod}(B,\mathbb{T}_{F,x}\otimes_B M) \simeq \mathbb{T}_{F,x}(M).$$
When $\mathbb{T}_{F,x}$ is moreover perfect, this implies that
$\mathbb{T}_{F,x}(M) \simeq Map_{B-Mod}(\mathbb{T}_{F,x}^{\vee},M),$
and thus that the cotangent complex of $F$ at $x$ exists and is
$\mathbb{L}_{F,x}=\mathbb{T}_{F,x}^{\vee}$.  \hfill $\Box$

\subsection{Infinitesimal properties of $\Perf^{\nabla}$}

We now study the infinitesimal structure of the derived moduli functors
$\Perf^{\nabla}(X)$ and $\Perf^{\nabla}(\fb X)$ constructed in the previous
section. The main result is the following. 

\begin{prop}\label{p5}
Let $X$ be a smooth algebraic variety over $k$ and let $X
\hookrightarrow \Xbar$ be a good compactification. Then, the two
derived moduli stacks $\Perf^{\nabla}(X)$ and $\Perf^{\nabla}(\fb X)$
are nilcomplete and infinitesimally cartesian.
\end{prop}

\

\noindent
\textit{Proof:} We start with $\Perf^{\nabla}(X)$. By construction
this derived stack is a derived mapping stack and can be written of
the form $\Perf^{\nabla}(X) \simeq \Map_{\dSt_k}(X_{DR},\Perf),$ where
$X_{DR}$ is the de Rham functor associated to $X$ (see for example
\cite{gr} for the relation between $\D$-modules and sheaves on
$X_{DR}$). We can write $X= \op{colim} Spec\, A_i$ as a finite colimit
of affine schemes, and thus $X_{DR}\simeq \op{colim} (Spec\,
A_i)_{DR}$.  The derived stack $\Perf^{\nabla}(X)$ is then the limit
of $\Perf^{\nabla}(Spec\, A_i)$. Since a limit of nilcomplete
(respectively infinitesimally cartesian) derived stacks is again
nilcomplete (respectively infinitesimally cartesian), we have reduced
the statement to the case where $X=Spec\, A$ is furthermore
afffine. The $\Perf^{\nabla}(X)$ statement thus boils down to the
following

\begin{lem}\label{l4}
Let $F$ be a nilcomplete (respectively infinitesimally cartesian)
derived stack over $k$. For any affine scheme $X$, the derived mapping
stack $\Map_{\dSt_k}(X,F)$ is again nilcomplete (respectively
infinitesimally cartesian).
\end{lem}

\

\noindent
\textit{Proof of the lemma:} Let $X=Spec\, A$, and $B$ any connective
cdga. Assume first that $F$ is nilcomplete. We consider the Postnikov
tower $\{B_{\leq n}\}_n$ of $B$.  The natural map
$$\Map_{\dSt_k}(X,F)(B) \longrightarrow \lim_n \Map_{\dSt_k}(X,F)(B_{\leq n})$$
can be written as 
$$F(A\otimes_k B) \longrightarrow \lim_n F(A\otimes_k B_{\leq n}).$$
As $k$ is a field, $A$ is a flat over $k$, and the tower $\{A\otimes_k
B_{\leq n}\}_n$ is a Postnikov tower for $A\otimes_k B$, and thus by
the assumption on $F$ the above morphism is an equivalence.

Let us now assume that $F$ is infintesimally cartesian. Let $B\oplus_d
M$ be a square zero extension of $B$ by a connective module $M$, given
by a cartesian square
$$\xymatrix{ B \oplus_d M \ar[r] \ar[d] & B \ar[d]^-{d} \\ B
  \ar[r]_-{0} & B\oplus M[1].}$$ Again because $A$ is flat over $k$,
tensoring with $A$ induces a pull-back diagram of connective cdga
$$\xymatrix{ C \oplus_d M_C \ar[r] \ar[d] & C \ar[d]^-{d} \\ C
  \ar[r]_-{0} & C\oplus M_C[1],}$$ where $C:=A\otimes_k B$ and
$M_C:=C\otimes_B M$. As $F$ is assumed infinitesimally cartesian, the
image of this diagram by $F$ remains a pull-back. By definition this
diagram is equivalent to
$$\xymatrix{ \Map_{\dSt_k}(X,F)(B\oplus_d M) \ar[r] \ar[d] &
  \Map_{\dSt_k}(X,F)(B) \ar[d] \\ \Map_{\dSt_k}(X,F)(B) \ar[r] &
  \Map_{\dSt_k}(X,F)(B\oplus M[1]).}$$ This shows that
$\Map_{\dSt_k}(X,F)$ is infintesimally cartesian.  \hfill $\Box$

\

\noindent
Next we analyse $\Perf^{\nabla}(\fb X)$. The argument here is slightly
different since this is not a derived mapping stack. We start by
writting $\Xbar= \op{colim} Spec\, A_i$ as a colimit of open affine
sub-schemes.  Without a loss of generality we can assume that the
divisor $D$ is principal on each $Spec\, A_i$, defined by an equation
$f_i \in A_i$. By the descent result of \cite{hpv}, we know that 
$\Perf^{\nabla}(\fb X)$ is then equivalent to a limit
$\Perf^{\nabla}(\fb X)=\lim_i F_i$ of derived stacks.  These derived
stacks $F_i$ can be described as follows. For each connective cdga
$B$, we have the completeted de Rham algebra of $A_i\otimes_k B$
defined by
$$\widehat{DR}_B(A_i):=\lim_j (DR(A_i/(f_i)^j) \otimes_k B).$$
This is a
$B$-linear graded mixed cdga for which the weight zero part is
$\widehat{A_i\otimes_k B}:=\lim_j (A_i/(f_i)^j \otimes_k B)$. By
inverting the weight zero element $f_i$ we have a new graded mixed
cdga
$$\widehat{DR}^{o}_B(A_i):=\lim_j (DR(A_i/(f_i)^j) \otimes_k
B)[f_i^{-1}].$$
The derived stack $F_i$ is then the functor sending
$B$ to the space of all graded mixed $\widehat{DR}^{o}_B(A_i)$-dg-modules
which are perfect as $\widehat{A\otimes_k B}$-dg-modules. Let us drop
the index $i$ and simply write $A$ and $f$ for  $A_i$ and 
$f_i$.  The derived stacks under consideration naturally carry
structures of stacks in dg-categories and will be considered as such
below.

We have a forgetful dg-functor
$$\widehat{DR}^{o}_B(A)-\medg \longrightarrow
\widehat{DR}^{o}_B(A)-\dg$$ from graded mixed dg-modules to
dg-modules. According to Corollary~\ref{cp2} this realizes the left
hand side as the dg-category of fixed points in
$\widehat{DR}^{o}_B(A)-\dg$ for the natural action by the group
$\mathcal{H}$ on the right hand side.  Restricting to dg-modules which
are perfect over $\widehat{A\otimes_k B}$ on both sides provides a
similar forgetful dg-functor
$$F(A) \longrightarrow Perf(\widehat{DR}^{o}_B(A))$$. Our statement
now reduces to the following

\begin{lem}\label{l5}
The $\s$-functor $B \mapsto \iPerf(\widehat{DR}^{o}_B(A))$ is
nilcomplete and infinitesimally cartesian as a derived stack of
dg-categories.
\end{lem}

\

\noindent
\textit{Proof of the lemma:} Note that if $\{B_{\leq
  n}\}_n$ is the Postnikov tower for $B$, then $\{\widehat{A\otimes_k
  B_{\leq n})}\}_n$ is a Postnikov tower for $\widehat{A\otimes_k
  B}$. In the same manner, $\widehat{A\otimes_k (-)}$ will transform a
square zero extension to a square zero extension.  The lemma therefore
reduces to the fact that the derived stack $\Perf$ is nilcomplete and
infinitesimally cartesian. This however is automatic since $\Perf$ is
locally geometric (see \cite{tv}).  \hfill $\Box$

\

\noindent
Finally the proof of the proposition follows from lemma \ref{l5} and
the fact that the operation of taking $\mathcal{H}$-fixed points 
preserves limits of dg-categories.  \hfill $\Box$

\medskip

\noindent
It is unclear to us whether $\Perf^{\nabla,ex}(\fb X)$ is also
nilcomplete and infinitesimally cartesian. Again, we beleive that the
inclusion $\Perf^{\nabla,ex}(\fb X) \hookrightarrow \Perf^{\nabla}(\fb
X)$ is an equivalence, but are unable to prove this at the
moment. Note also that we explicitly included the compactification
$\Xbar$ in the statement of the proposition above as we do not know
that $\Perf^{\nabla}(\fb X)$ is independent of the choice of
compactification (as opposed to $\Perf^{\nabla,ex}(\fb X)$).

\subsection{Cotangent complexes}

We now turn to the study of cotangent complexes of the derived stacks
$\Perf^{\nabla}(X)$ and $\Perf^{\nabla}(\fb X)$. In general, these
cotangent complexes do not exist, except when $X$ is proper. We thus
introduce the following notion.

\begin{df}\label{d11}
Let $B$ be a connective cdga. We say that an object $E \in
\Perf^{\nabla}(X)(B)$ (respectively $E \in \Perf^{\nabla}(\fb X)(B)$)
is \emph{\bfseries End-Fredholm} (or simply \emph{\bfseries Fredholm})
if the cotangent complex at $E$ exists and is perfect.
\end{df}

\

\noindent
The computation of the tangent complexes of $\Perf^{\nabla}(X)$ and
$\Perf^{\nabla,ex}(\fb X)$ is standard and is given by the following
proposition.

\begin{prop}\label{pd11}
Let $B$ be a connective cdga and $E \in \Perf^{\nabla}(X)(B)$ with
restriction \linebreak $\mathsf{R}(E) \in \Perf^{\nabla,ex}(\fb X)(B)$.
\begin{enumerate}

\item The $\s$-functor $M \mapsto \mathbb{T}_{\Perf^{\nabla}(X),E}(M)$
  is equivalent to $M \mapsto \HH_{DR}(X,E^\vee\otimes E \otimes_B
  M)[1]$.

\item The $\s$-functor $M \mapsto \mathbb{T}_{\Perf^{\nabla}(\fb
  X),\mathsf{R}(E)}(M)$ is equivalent to $M \mapsto \HH_{DR}(\fb
  X,\mathsf{R}(E^\vee\otimes E) \otimes_B M)[1]$.

\end{enumerate}

\end{prop}

\

\noindent
As a direct consequence of the above proposition and the definition of
being Fredholm we have the following direct corollary.

\begin{cor}\label{cd11}
Let $B$ be any connective cdga and $E \in \Perf^{\nabla}(X)(B)$. 
\begin{enumerate}

\item[(1)] The object $E$ is Fredholm if and only if 
$\HH_{DR}(X,E^\vee \otimes E)$ is a perfect $B$-module.

\item[(2)] The object $\mathsf{R}(E)$ is Fredholm if and only if 
$\HH_{DR}(\fb X,\mathsf{R}(E^\vee \otimes E))$ is a perfect $B$-module.

\end{enumerate}

\end{cor}

\

\noindent
\textit{Proof:} $(1)$ By definition it is enough to show that the
formation of $\HH_{DR}(X,E^\vee \otimes E)$ is compatible with base
changes of $B$. But this follows immediately from the corresponding
statement for Tate $B$-module that was proven in
Section~\ref{s-cptsup}. The proof of $(2)$ is similar. \hfill $\Box$

\

\noindent
We will see later that all objects are Fredholm when $B$ is a field
(see corollary \ref{c7}). More generally, if $E \in
\Perf^{\nabla}(X)(B)$ is any object, we will see that $E$ and $R(E)$
are Fredholm under the condition that fro some good compactification
$j : X \hookrightarrow \Xbar$ we have that both $j_*(E)$ and
$j_*(E^\vee)$ are perfect $\D_X\otimes_k B$-modules in the sense of
Section~\ref{app:relD}.  We refer to corollary \ref{c6} for this
important statement which will be crucial in the proof of the
representability theorem.

\section{The Lagrangian restriction map}

In  this section we construct a natural shifted Lagrangian structure on the
restriction morphism
$$
\mathsf{R} : \Perf^{\nabla}(X) \longrightarrow \Perf^{\nabla}(\fb X),
$$
which is the de Rham analogue of the   Betti
statements in \cite{ptbetti}. However, the new feature
here is that the derived stacks $\Perf^{\nabla}(X)$ and
$\Perf^{\nabla}(\fb X)$ are not representable, and their tangent
complexes can be infinite dimensional. We thus have to be careful with
the notion of Lagrangian structure itself. The definitions of closed
forms and isotropic structures make sense on
general derived stacks. However, the non-degeneracy condition in the
definition of a Lagrangian structure causes a problem as there is a
priori no direct relationship between 2-forms on a derived stack $F$ and
global sections of $\wedge^2 \mathbb{L}_F$ (even assuming that
$\mathbb{L}_F$ exists).

Therefore in our setting non-degeneracy has to be defined pointwise, at
all field valued points. For this we use in a crucial manner that the
derived stacks $\Perf^{\nabla}(X)$ and $\Perf^{\nabla}(\fb X)$ are
nil-complete and infinitesimally cartesian, and moreover that their
cotangent complexes exist and are perfect at all field valued points
(see proposition \ref{p5} corollary \ref{c7})

\subsection{Closed forms and symplectic structures}

Recall from \cite{ptvv} that for any derived stack $F$ we have a
complex of $p$-forms $\A^p(F)$, and a complex of closed $p$-forms
$\A^{p,cl}(F)$, together with a forgetful morphism $\A^{p,cl}(F)
\longrightarrow \A^p(F)$. When $F=Spec\, A$ is a derived affine
scheme, the complex $\A^p(F)\sim \wedge_A^p \mathbb{L}_A$ simply is
the $p$-th wedge power of the cotangent complex of $A$. In the same
manner, $\A^{p,cl}(F)\sim \text{tot}(\prod_{i\geq
  p}(\wedge^i_A\mathbb{L}_A)[-i])$ is the totalization of the
completed derived truncated de Rham complex.

Suppose that $F$ is any derived stack that possesses a cotangent
complex $\mathbb{L}_F \in \Dqcoh(F)$ in the sense recalled in the
definition \ref{d9}.  There is a descent morphism
$$\mathbb{H}(F,\wedge^p_{\OO_F}\mathbb{L}_F) \longrightarrow \A^p(F).$$
When $F$ is a derived Artin stack, it is shown in \cite[Proposition~1.14]{ptvv} that
this morphism is a quasi-isomorphism. In general this descent morphism
has no reason to be a quasi-isomorphism. This fact creates
complications when one tries to define the non-degeneracy 
of $2$-forms \cite{ptvv}. In this paper we overcome this complication
by working pointwise on $F$ as follows.

\begin{df}\label{d12}
  A derived stack $F$ is \emph{formally good} if it is infinitesimally
  cartesian and for any $k$-field $L$ and any $x\in F(L)$ the tangent
  complex $\mathbb{T}_xF$ is perfect over $L$.
\end{df}

\

\noindent
Proposition \ref{p5} and corollary \ref{c7} show that 
the derived stacks $\Perf^{\nabla}(X)$ and $\Perf^{\nabla}(\fb X)$
are formally good in the sense of Definition~\ref{d12}. 

Let $F$ be a formally good derived stack and $x\in F(L)$ be a field
valued point.  We can restrict the functor $F$ to the $\s$-category of
artinian local augmented $L$-cdga by sending such a cdga
$A \in \dgart^*_L$ to the fiber of $F(A) \longrightarrow F(L)$ taken
at $x$.  By definition this restriction is the formal completion of
$F$ at $x$ and we will denote it by $\widehat{F}_{x}$.  Since $F$ is
assumed to be infinitesimally cartesian the $\s$-functor
$\widehat{F}_x$ is a formal moduli problem over $L$ in the sense of
\cite{lu3}.  It therefore corresponds to an $L$-linear dg-Lie algebra
$\LL_x$ whose underlying complex is $\mathbb{T}_xF[-1]$.

By left Kan extension from artinian cdga to connective cdga, the
$\s$-functor $\widehat{F}_x$ can be itself considered as a derived
stack. As such it possesses a complex of $p$-forms
$\A^p(\widehat{F}_x)$. It turns out that this complex can be computed
purely in terms of the dg-Lie algebra $\LL_x$ as follows.

\begin{prop}\label{p7}
  Let $F$ be a formal moduli problem over $L$, associated to a dg-Lie
  algebra $\LL$.  There is a canonical quasi-isomorphism
$$\A^p(F)\simeq \underline{Hom}_{\LL-\dg}(k,\wedge^p(\LL^\vee[-1])),$$
where $\LL^{\vee}$ is the $L$-linear dual of $\LL$ considered as a
dg-module over $\LL$ by the coadjoint action.
\end{prop}

\

\noindent
\textit{Proof:} We first prove the statement when $F$ is
representable, that is $F=Spec\, A$ for $A \in \dgart^*_L$. In this
case $F$ has a cotangent complex $\mathbb{L}_{A/L} \in \Dqcoh(F)$. By
\cite{lu3} there is a full embedding \linebreak
$\Dqcoh(F) \hookrightarrow D(\LL-\dg)$ and the image of
$\mathbb{L}_{A/L}$ is the dg-module $\LL^{\vee}[-1]$, which follows
immediately from the universal property of
$\mathbb{L}_{A/L}$. Finally, the above full embedding also sends $\OO$
to $k$, which implies the existence of the required equivalence
$$\A^p(F)=\underline{Hom}(\OO,\wedge^p\mathbb{L}_{A/L})
\simeq \underline{Hom}_{\LL-\dg}(k,\wedge^p(\LL^\vee[-1])).$$ This
extends easily to the case where $F=\op{colim} Spec\, A_i$ is now only
pro-representable by a pro-object $"\lim_i A_i"$ in $\dgart^*_L$.

To deduce the general case we use the existence of smooth
hyper-coverings proved in \cite{lu3}. Having smooth hyper-coverings
guarantees that a general formal moduli problem $F$ can be written as
a geometric realization $|F_*|$ of a simplicial object in
pro-representables which moreover satisfies the smooth hyper-coverings
condition. We can then use the same descent argument as done in the
algebraic case in \cite{ptvv}.  We consider the formal moduli problem
$TF[-1]=Map(Spec\, (k\oplus k[1]),F)$ corresponding to the shifted
tangent of $F$. Using a smooth hyper-covering $F_{*}$ we observe that
$TF[-1]$ is again the realization of $TF_*[-1]$. Passing to the
complex of functions we find that the natural morphism
$$\underline{Hom}_{\LL-\dg}(k,\wedge^p(\LL^\vee[-1])) \longrightarrow
\lim_n \underline{Hom}_{\LL_n-\dg}(k,\wedge^p(\LL_n^\vee[-1]))$$ is a
quasi-isomorphism (where we have denoted by $\LL_n$ the dg-Lie algebra
corresponding to $F_n$). This last descent statement, together with
the already treated case of pro-representable $F$ proves the general
result. \hfill $\Box$

\

\noindent
Going back to our formally good stack $F$ let $x \in F(L)$ be a field
valued point. Using proposition \ref{p7} we see that there is a
natural restriction map
$$\A^p(F) \longrightarrow \A^p(\widehat{F}_x) 
\simeq \underline{Hom}_{\LL_x-\dg}(k,\wedge^p(\LL_x^\vee[-1]))
\longrightarrow \wedge^p(\LL_x^\vee[-1]),$$ where the last morphism is
obtained by forgetting the $\LL_x$-module structure.

\begin{df}\label{d13}
  Let $F$ be a formally good derived stack and
  $\omega \in H^n(\A^{2,cl}(F))$ be a closed $2$-form of degree $n$ on
  $F$. We say that \emph{$\omega$ is non-degenerate} if for all field
  valued points $x\in F(L)$, the image of $\omega$ by the morphism
  $$\A^{2,cl}(F) \longrightarrow A^{2}(F)
  \longrightarrow \wedge^p(\LL_x^\vee[-1])\simeq 
\wedge^2(\mathbb{T}^{\vee}_{F,x})$$ is a non-degenerate pairing of
degree $n$ and induces an equivalence
$\mathbb{T}_{F,x} \simeq \mathbb{T}_{F,x}^{\vee}[n].$
\end{df}

\

\noindent
The above definition generalizes immediately to the relative setting as
follows.  Suppose now that we have a morphism of formally good derived
stacks $f : F \longrightarrow F'$, and $\omega$ a closed $2$-form of
degree $n$ on $F'$. Assume that we are given a homotopy to zero
$h : f^*(\omega) \sim 0$ inside $\A^{2,cl}(F)$. By what we have seen,
for any field valued point $x \in F(L)$, the form $\omega$ and the
homotopy $h$ induces an $n$-cocycle $\omega_x$ in
$\wedge^2\mathbb{T}^{\vee}_{F',x}$ as well as a null homotopy of its
image $f^*(\omega_x)$ in $\wedge^2\mathbb{T}_{F,x}^{\vee}$.  This null
homotopy induces a well defined morphism of complexes
$$\mathbb{T}_{F,x} \longrightarrow \mathbb{T}_{F/F',x}^{\vee}[n-1].$$

\begin{df}\label{d14}
  Let $f : F \longrightarrow F'$ be a morphism of formally good
  derived stacks.  Let $\omega$ be a closed $2$-form of degree $n$ on
  $F'$ and $h: f^*(\omega) \sim 0$ an isotropic structure on $f$
  with respect to $\omega$. We say that the isotropic structure is
  \emph{Lagrangian} if for any field valued point $x\in F(L)$ the
  induced morphism of complexes
$$\mathbb{T}_{F,x} \longrightarrow \mathbb{T}_{F/F',x}^{\vee}[n-1]$$
is a quasi-isomorphism.
\end{df}

\subsection{Orientation on the formal boundary} \label{s-orientation}

In this section we will prove that the conditions for applying the
results of \cite{to3} are satisfied for the restriction morphism of
derived stacks
$$
\mathsf{R} : \Perf^{\nabla}(X) \longrightarrow \Perf^{\nabla}(\fb X).
$$
The main step consists of studying Serre duality on $\fb X$ and the
key ingredient is the construction of the integration map
$$
or : \HH(\fb X,R(\omega_X)) \longrightarrow k[1-d]
$$
where $d$ is the dimension of $X$ (for simplicity we assume that
$X$ is connected).

Here $\omega_X:=\Omega_X^d$ is the canonical sheaf of $X$. We pick a
good compactification $j : X \hookrightarrow \Xbar$ once and for
all. As before we will write $\widehat{\Xbar}$ for the formal
completion of $\Xbar$ along the divisor $D = \Xbar - X$ and we will
write $\widehat{\jfrak} : \widehat{\Xbar} \to \Xbar$ for the natural map.

The formal gluing theorem of \cite{hpv} and the observation that
$\omega_{X} = j^{*}\omega_{\Xbar}$, yield a cartesian sqare
\begin{equation} \label{eq:glue.square}
  \xymatrix{
  \HH(\Xbar,\omega_{\Xbar}) \ar[d] \ar[r] & \HH(X,\omega_X) \ar[d] \\
  \HH(\widehat{\Xbar},\widehat{\jfrak}^{*}{\omega}_{\Xbar}) \ar[r] & \HH(\fb
  X,R(\omega_X)).}
\end{equation}
The boundary map for this cartesian square produces a morphism $u :
\HH(\fb X,R(\omega_X)) \longrightarrow \HH(\Xbar,\omega_{\Xbar})[1]$ of
complexes over $k$.  Composing with Grothendieck's trace isomorphism
$H^d(\Xbar,\omega_{\Xbar})\simeq k$ we get the required morphism of
complexes
$$or : \HH(\fb X,R(\omega_X)) \longrightarrow k[1-d].$$
This morphism is a version of the residue map, for instance it
coincides with the usual residue of forms when $X$ is a curve and the
residues are taken at the points at infinity.

We defined the morphism $or$ above as a morphism of complexes over
$k$. However, as explained in Section~\ref{s-cptsup}, the source of
this morphism is the realization of the ind-pro complex
$\widetilde{\HH}(\fb X,R(\omega_X))$. By construction the formal gluing
giving the cartesian square \eqref{eq:glue.square} lifts canonically
to give a cartesian square of Tate complexes over $k$. This implies
that the boundary morphism $u$ also lifts canonically
as a morphism in the ind-pro category.  As a result $or$ arises as the
realization of a natural moprhism of Tate complexes
$$
\widetilde{or} : \widetilde{\HH}(\fb X,R(\omega_X)) \longrightarrow k[1-d].
$$
By base change (see Section~\ref{s-cptsup}) we get an induced
morphism for every connective cdga $B$
$$
\widetilde{or} : \widetilde{\HH}(\fb X,R(\omega_X)\otimes_k B)
\longrightarrow B[1-d].
$$
Assume now that $B$ is a connective cdga and $E$ and $F$ are two
perfect complexes over $X\times S$, with $S=Spec\, B$. To simplify the
discussion we assume that $E$ and $F$ can be extended to perfect
complexes on $\Xbar\times S$ (even though this is not strictly
necessary for the results below).  We have a composition morphism
$$
\HH(\fb X,R(E)^\vee \otimes R(F)) \otimes_B \HH(\fb X,R(F)^\vee
\otimes R(E)\otimes R(\omega_X)) \longrightarrow \HH(\fb X,R(E)^\vee
\otimes R(E)\otimes R(\omega_X))$$ which we can compose with the trace
morphism $R(E)^\vee \otimes R(E) \longrightarrow R(\OO_X)$, and with
the orientation $or$ in order to get a pairing
$$
\HH(\fb X,R(E)^\vee \otimes R(F)) \otimes_B \HH(\fb X,R(F)^\vee
\otimes R(E)\otimes R(\omega_X)) \longrightarrow B[1-d].
$$
This pairing also admits a canonical lift as a pairing of Tate
$B$-modules. Indeed, we already have seen that $or$ has such a lift,
and composition and trace are also compatible with the ind-pro
structures. We thus have defined a canonical pairing of Tate
$B$-modules
$$
\widetilde{\HH}(\fb X,R(E)^\vee \otimes R(F)) \widehat{\otimes}_B
\widetilde{\HH}(\fb X,R(F)^\vee \otimes R(E)\otimes R(\omega_X))
\longrightarrow B[1-d].
$$
By rigidity we may assume that $F=\OO_X$ without loss of
generality. The pairing can then be written as
$$
\widetilde{\HH}(\fb X,R(E)^\vee) \widehat{\otimes}_B
\widetilde{\HH}(\fb X,R(E)\otimes R(\omega_X)) \longrightarrow B[1-d].
$$
By construction, the orientation morphism $\widetilde{or}$ canonically
vanishes on $\HH(X,\omega_X)$, and so we get 
an induced pairing of Tate $B$-modules
\begin{equation} \label{eq:serre.pairing}
\HH(X,E^\vee) \widehat{\otimes}_B \HHc(X,E\otimes \omega_X)
\longrightarrow B[-d].
\end{equation}
The following result is Serre duality for cohomology with compact
supports.

\begin{prop}\label{p8}
  The pairing \eqref{eq:serre.pairing}
  is non-degenerate. It induces an equivalence of Tate $B$-modules
  $$
  \HHc(X,E\otimes \omega_X) \simeq \HH(X,E^\vee)^{\vee}[-d].
  $$
\end{prop}

\

\noindent
\textit{Proof:} The pairing \eqref{eq:serre.pairing} induces a
morphism of $B$-modules $\alpha : \HH(X,E^\vee) \longrightarrow
\HHc(X,E\otimes \omega_X)^\vee[-d].$ Here, $\HHc(X,E\otimes
\omega_X)^\vee$ is the dual of $\HHc(X,E\otimes \omega_X)$ as a Tate
module. Since $\HHc(X,E\otimes \omega_X)$ is pro-perfect this dual is
a genuine $B$-module. We must show that the morphism $\alpha$ is an
equivalence. For this, we go back to examine the formal gluing
cartesian square \eqref{eq:glue.square}, and the definition of the
pairing. We have the exact triangle of Tate $B$-modules
$$
\xymatrix{\HHc(X,E\otimes \omega_X) \ar[r] & 
\HH(\Xbar,E\otimes \omega_{\Xbar}) \ar[r] &
\text{``$\lim_n$"} \, \HH(\Xbar_{(n)},\jfrak_n^* (E\otimes \omega_{\Xbar})).}
$$
The rightmost term  can be written as 
$\text{``$\lim_n$"}\, \HH(\Xbar_{(n)},\jfrak_n^*(E)\otimes \omega_{\Xbar_{(n)}}
\otimes \LL_{n})$, where
$\LL_{n}$ is the conormal sheaf of $\jfrak_{n} : \Xbar_{(n)} \hookrightarrow
\Xbar$.
By Serre duality on $\Xbar$ and $\Xbar_{(n)}$ (for each $n$),
the restriction map 
$$
\HH(\Xbar,E\otimes \omega_X) \longrightarrow
\HH(\Xbar_{(n)},\jfrak_n^*(E)\otimes \omega_{\Xbar_{(n)}} \otimes \LL_{n})
$$
is dual to the natural map
$\HH(\Xbar_{(n)},\jfrak_n^*(E^{\vee})\otimes \LL_{n}^{\vee})[d-1]
\longrightarrow \HH(\Xbar,E^{\vee})[d].$
Passing to the colimit over $n$ these assemble in a  natural map
$$\HH_{D}(\Xbar,E^{\vee})[d] \longrightarrow \HH(\Xbar,E^{\vee})[d],$$
where the source is cohomology with supports in $D$. The cofiber of
this map is then naturally equivalent to $\HH(X,E^{\vee})[d]$. This
constructs a natural equivalence of $B$-modules
$$
\HHc(X,E\otimes \omega_X)^{\vee} \simeq \HH(X,E^\vee)[d].
$$ It is straightforward to check that this equivalence is the
morphism $\alpha$.  \hfill $\Box$

\

\begin{cor}\label{c5}
The pairing of Tate $B$-modules 
$$
\widetilde{\HH}(\fb X,R(E)^\vee) \widehat{\otimes}_B
\widetilde{\HH}(\fb X,R(E)\otimes R(\omega_X))
\longrightarrow B[1-d]
$$
is non-degenerate.
\end{cor}

\

\noindent
\textit{Proof:} We have two exact triangles of Tate $B$-modules
$$
\xymatrix{
  \HHc(X,E^\vee) \ar[r] & \HH(X,E^{\vee}) \ar[r] &
  \widetilde{\HH}(\fb X,R(E^{\vee}))}
$$
and
$$
\xymatrix{ \HHc(X,E\otimes \omega_X) \ar[r] & \HH(X,E\otimes \omega_X)
  \ar[r] & \widetilde{\HH}(\fb X,R(E \otimes \omega_X))}
$$
The dual,
inside Tate $B$-modules, of the second triangle is (up to a rotation
and shift by $-d$)
$$
\xymatrix{\HH(X,E\otimes \omega_X)^{\vee}[-d] \ar[r] & \HHc(X,E\otimes
  \omega_X)^{\vee}[-d]\ar[r] & \widetilde{\HH}(\fb X,R(E\otimes
  \omega_X))^{\vee} [1-d]}.
$$
By the construction of the orientation $or$ the natural pairing
produces a commutative diagram of Tate $B$-modules
$$
\xymatrix{ \HHc(X,E^\vee) \ar[d] \ar[r] & \HH(X,E^{\vee}) \ar[d]
  \ar[r] & \widetilde{\HH}(\fb X,R(E^{\vee})) \ar[d] \\ \HH(X,E\otimes
  \omega_X)^{\vee}[d] \ar[r] & \HHc(X,E\otimes
  \omega_X)^{\vee}[d]\ar[r] & \widetilde{\HH}(\fb X,R(E\otimes
  \omega_X))^{\vee}[d-1].}
$$
The first two vertical morphisms on the left are equivalences by
proposition \ref{p8}.  Therefore the third vertical morphism is also
an equivalence. \hfill $\Box$

\

\medskip

\noindent
The same orientation morphism can be used to prove a duality statement
for de Rham cohomology with compact supports. It goes as follows. The
complex of sheaves $|\widehat{\DR}_{B}^{o}|$ on $D$,
computing\footnote{Here we use a slight abuse of notation and write
  simply $\OO_{X}$ for the trivial rank one flat bundle
$(\OO_{X},d_{DR})$ on $X$.} $\HH_{DR}(\fb X,\mathsf{R}(\OO_X))$ is
  bounded of amplitude contained in $[0,d]$. Moreover, its last
  non-zero term is $R(\omega_X)$. Therefore, there is a canonical map
$$
  H^{2d-1}_{DR}(\fb X,\mathsf{R}(\OO_X)) \longrightarrow
  H^{d-1}(\fb X,R(\omega_X)).
$$
Composing with the orientation map $or : H^{d-1}(\fb X,R(\omega_X))
\longrightarrow k[1-d]$ we get an orientation morphism $\HH_{DR}(\fb
X,\mathsf{R}(\OO_X)) \longrightarrow k[1-2d]$. As before it extends
naturally
as a moprhism of Tate complexes \linebreak over $k$
$$
\widetilde{or} : \widetilde{\HH}_{DR}(\fb X,\mathsf{R}(\OO_X))
\longrightarrow k[1-2d].
$$
For any connective cdga $B$ and any $E \in \Perf^{\nabla}(X)(B)$, this
orientation defines as before two pairings of Tate $B$-modules
\begin{equation} \label{eq:pairing1}
\HHcd(X,E) \widehat{\otimes}_B \widetilde{\HH}_{DR}(X,E^\vee)
\longrightarrow B[-2d]
\end{equation}
and 
\begin{equation} \label{eq:pairing2}
\widetilde{\HH}_{DR}(\fb X,\mathsf{R}(E)) \widehat{\otimes}_B
\widetilde{\HH}_{DR}(\fb X,\mathsf{R}(E)^\vee) \longrightarrow B[1-2d].
\end{equation}
We now have the following

\begin{prop}\label{p9}
  The pairings \eqref{eq:pairing1} and \eqref{eq:pairing2}
  are non-degenerate and induce natural equivalences
of Tate $B$-modules
$$
\HHcd(X,E) \simeq \widetilde{\HH}_{DR}(X,E^\vee)^\vee[1-2d] \qquad 
\widetilde{\HH}_{DR}(\fb X,\mathsf{R}(E)) \simeq 
\widetilde{\HH}_{DR}(\fb X,\mathsf{R}(E)^\vee)^\vee[-2d].
$$
\end{prop}

\

\noindent
\textit{Proof:} We use the Hodge filtrations on the
various complexes computing these cohomology groups. In terms of
graded mixed modules these are the filtrations on $|E|$ given by
$\oplus_{i\geq p}E(i)[-2i] \subset \oplus_{i}E(i)[-2i]$. The
associated graded of these filtrations are perfect complexes of the
form $E(0) \otimes_{\OO_X}\Omega_X^i[-i]$. The pairings of the
proposition are compatible with these filtrations and the induced
pairings are the one for Serre duality of perfect complexes. Therefore
the proposition follows from the  Serre duality with compact
supports from Proposition \ref{p8}. \hfill $\Box$

\

\medskip

\noindent
One important corollary of the previous results is the following 
criterion for finiteness of the de Rham cohomology of $\fb X$. 

\begin{cor}\label{c6}
Let $E \in \Perf^{\nabla}(X)(B)$ be such that $\HH_{DR}(X,E)$ and
$\HH_{DR}(X,E^{\vee})$ are both perfect \linebreak $B$-modules. Then
the Tate $B$-modules $\widetilde{\HH}_{DR}(\fb(X),\mathsf{R}(E))$
and $\HHcd(X,E)$
are both perfect.
\end{cor}

\

\noindent
\textit{Proof:} Using the exact triangle
$$\xymatrix{ \HHcd(X,E) \ar[r] & \HH_{DR}(X,E) \ar[r] &
  \widetilde{\HH}_{DR}(\fb X,\mathsf{R}(E))}$$ we see that the Tate $B$-module
$\widetilde{\HH}_{DR}(\fb X,\mathsf{R}(E))$ must be pro-perfect. But corollary
\ref{c5} implies that its dual is also pro-perfect. This implies that
it must be perfect. \hfill $\Box$

\

\noindent
One important consequence is the following. 

\

\begin{cor}\label{c7}
\begin{enumerate}
\item[(1)]
Let $E \in \Perf^{\nabla}(X)(B)$ be such that $\HH_{DR}(X,E^\vee\otimes E)$ 
is perfect over $B$. Then $E$ and $\mathsf{R}(E)$ are both Fredholm in the sense of 
definition \ref{d11}.
\item[(2)] If $B=k$, any $E \in \Perf^{\nabla}(X)(k)$ is Fredholm and
  so is $\mathsf{R}(E)$.
\end{enumerate}
\end{cor}

\

\noindent
\textit{Proof:} $(1)$ is a direct consequence of corollary \ref{c6}
and the fact that both $\HH_{DR}(X,E)$ and
$\widetilde{\HH}_{DR}(\fb(X),R(E))$ are stable by base changes of $B$.  For
$(2)$, we have to show that for any $E\in \Perf^{\nabla}(X)(k)$ the
complex $\HH_{DR}(X,E^{\vee}\otimes E)$ is perfect over $k$.  But
$E^\vee \otimes E$ is a bounded complex of coherent $\D_X$-modules
with holonomic cohomologies. By Bernstein's theorem holonomic
$\D$-modules are stable by push-forward and so $\HH_{DR}(X,E^{\vee}
\otimes E)$ is a bounded complex with finite dimensional cohomology
and thus perfect. \hfill $\Box$

\

\noindent
We are now are ready to construct the Lagrangian structure on
the restriction morphism
$$
\mathsf{R} : \Perf^{\nabla}(X) \longrightarrow
\Perf^{\nabla}(\fb X).
$$
For this we use the main result of
\cite{to3}. The derived stack $\Perf^{\nabla}(X)$ is  the
underlying stack of a derived stack in symmetric monoidal rigid
dg-categories. According to \cite{to3} in order to construct a closed
$2$-form $\omega$ on $\Perf^{\nabla}(\fb X)$, together with an
homotopy $h : \mathsf{R}^*(\omega) \sim 0$, it is enough to:

\begin{enumerate}
\item[(i)] construct a morphism of complexes of $k$-modules
$$or : \HH_{DR}(\fb X,\mathsf{R}(\OO_X)) \longrightarrow k[1-2d]$$
together with a homotopy to zero of the restriction
$$\mathsf{R}^*(or) : \HH_{DR}(X,\OO_X) \longrightarrow k[1-2d]$$
and 
\item[(ii)] prove that for any connective cdga $B$ the induced morphisms
  $$
\HH_{DR}(\fb X,\mathsf{R}(\OO_X))\otimes_k B \longrightarrow \HH_{DR}(\fb
X,\mathsf{R}(\OO_X)\otimes_k B)
$$
and
$$\HH_{DR}(X,\OO_X)\otimes_k B \longrightarrow \HH_{DR}(X,\OO_X\otimes_k B)$$
are equivalences of $B$-modules.
\end{enumerate}

\

\noindent
Statement $(ii)$ holds thanks to corollary \ref{c6}. The map $or$ is
constructed at the beginning of the section. Recall that it
comes from the cartesian square
$$
\xymatrix{
\HH(\Xbar,\omega_{\Xbar}) \ar[r] \ar[d] & \HH(X,\omega_X) \ar[d] \\
\HH(\widehat{\Xbar},\widehat{\jfrak}^{*}\omega_{\Xbar}) \ar[r] &
\HH(\fb X,R(\omega_X)),}
$$
and the associated  boundary map
$H^{d-1}(\fb X,R(\omega_X)) \longrightarrow
H^{d}(\Xbar,\omega_{\Xbar})\simeq k$. Precomposing with 
the canonical map $H^{2d-1}_{DR}(\fb X,\mathsf{R}(\OO_X)) \longrightarrow
H^{d-1}(\fb X,R(\omega_X))$ provides the orientation morphism
$$or : \HH_{DR}(\fb X,\mathsf{R}(\OO_X)) \longrightarrow k[1-2d].$$ By
construction, the composition $H^{d-1}(X,\omega_X) \longrightarrow
H^{d-1}(\fb X,R(\omega_X)) \longrightarrow
H^{d}(\Xbar,\omega_{\Xbar})$ is the zero map so a null homotopy of the
morphism $\mathsf{R}^*(or) : \HH_{DR}(X,\OO_X) \longrightarrow k[1-2d]$ is
given by a morphism $H^{2d}_{DR}(X,\OO_X) \longrightarrow k$. If $X$
is proper, we take this map to be the natural isomorphism. If $X$ is
not proper then $H^{2d}_{DR}(X,\OO_X)=0$ and this map is the zero map.

By the main result of \cite{to3}, we have that the derived stack
$\Perf^{\nabla}(\fb X)$ carries a canonical closed $2$-form $\omega$
of degree $3-2d$. Moreover the pull-back form $\mathsf{R}^*(\omega)$
comes equiped with a natural null-homotopy $h :
\mathsf{R}^*(\omega)\sim 0$. We thus have proved the following
statement.

\begin{cor}\label{c8}
The morphism of derived stacks $\mathsf{R} : \Perf^{\nabla}(X) \longrightarrow
\Perf^{\nabla}(\fb X)$ carries a canonical isotropic structure of
degree $2-2d$.
\end{cor}

As explained in definition \ref{d14}, the non-degeneracy condition on
an isotropic structure is imposed at all field valued points of
$\Perf^{\nabla}(X)$. Given such point $E \in \Perf^{\nabla}(X)(L)$
defined over a $k$-field $L$, the morphism 
$$\mathbb{T}_{\Perf^{\nabla}(X),E} \longrightarrow
\mathbb{L}_{\Perf^{\nabla}(X)/\Perf^{\nabla}(\fb X),E}[2-2d]
$$
induced
by the isotropic structure becomes, after the identifications given by
proposition \ref{pd11}, equal to the duality morphism
$$
\HH_{DR}(X,E^\vee\otimes E) \longrightarrow
\HHcd(X,E^\vee\otimes E)^{\vee}[-2d].
$$
The latter morphism is an equivalence by proposition \ref{p9}. This
proves the following

\begin{cor}\label{c9}
The isotropic structure of corollary \ref{c8} is a Lagrangian
structure in the sense of definition \ref{d14}.
\end{cor}

\

\section{The relative representability theorem} \label{sec-represent}

In this section we prove that the fibers of the restriction morphism
$\mathsf{R}$ over field valued points are locally representable by
quasi-algebraic spaces in the sense of our
Appendix~\ref{sec:artin-lurie}.  We prove this statement for vector
bundles endowed with flat connections.  The extension to the perfect
complexes setting can be reduced to this special case by truncation
and we leave it to the interested reader to fill in the details.

We first consider the derived substack $\Vect^{\nabla}(X) \subset
\Perf^{\nabla}(X)$ consisting of all objects whose underlying
$\OO_X$-module is a vector bundle. Explicitly, for a connective cdga
$B$, an object $E \in \Perf^{\nabla}(X)(B)$ lies in $\Vect^{\nabla}(X)
(B)$ if the $\OO_X\otimes_k B$-module $E(0)$ is locally free of finite
rank. We define similarly $\Vect^{\nabla}(\fb X)(B) \subset
\Perf^{\nabla}(\fb X)(B)$ as objects $E$ such that $E(0)$ is locally
free of finite rank as a $\widehat{\OO}_{D,B}^{o}$-module.

We fix once for all $V_\s \in \Vect^{\nabla}(\fb X)(k)$, a vector
bundle with flat connecion on the formal boundary of $X$. The fiber of
the restriction morphism $\mathsf{R} : \Vect^{\nabla}(X) \longrightarrow
\Vect^{\nabla}(\fb X)$ taken at $V_\s$ will be denoted by
$\Vect^{\nabla}_{V_{\s}}(X)$. It is the derived stack of vector
bundles with flat connecions on $X$ framed by $V_{\s}$ along $\fb X$.
When no component of $X$ is proper, the rank of $V_\s$ fixes the rank
of all objects in $\Vect^{\nabla}_{V_\s}(X)$.  Since the proper case
of the result is well understood we will assume that $X$ has no proper
component.

\begin{thm}\label{t1}
With the notations above, the derived stack $\Vect^{\nabla}_{V_\s}(X)$ is a
derived quasi-algebraic space in the sense of definition \ref{d17}.
\end{thm}

\

\noindent
\textit{Proof:} We will prove the theorem by applying the version of
Artin-Lurie representability criterion by quasi-algebraic derived
spaces recalled in Theorem~\ref{t2} of our
Appendix~\ref{sec:artin-lurie}. By Galois descent we may assume that
$k$ is algebraically closed. We also assume that the derived stack
$\Vect^{\nabla}_{V_\s}(X)$ is not empty, or equivalently that $V_\s$
extends to a flat vector bundle $V$ on the whole $X$.
 
By proposition \ref{p5} we know that $\Vect^{\nabla}_{V_\s}(X)$ is
infinitesimally cartesian and nil-complete, since it is defined as the
fiber of a morphism between two infinitesimally cartesian and
nil-complete derived stacks.  Let us show moreover that it has a
global cotangent complex. By Definition~\ref{d11} this amounts to
the following lemma.

\begin{lem}\label{l6}
Let $B$ be any connective cdga and $E\in \Vect^{\nabla}_{V_\s}(X)(B)$
an object. Then the image of $E$ in $\Vect^{\nabla}(X)(B)$ is Fredholm
over $B$.
\end{lem}

\

\noindent
\textit{Proof of the lemma:} This is a consequence of our corollary
\ref{cd11}. Indeed, we have an exact triangle of Tate $B$-modules
$$
\xymatrix{
\HHcd(X,E^\vee \otimes E) \ar[r] & \HH_{DR}(X,E^\vee \otimes E) \ar[r] & 
\widetilde{\HH}_{DR}(\fb X,R(E^\vee \otimes E)).}
$$ The rightmost module is equivalent to $\widetilde{\HH}_{DR}(\fb
X,R(V_\s^\vee \otimes V_\s))\otimes_k B$ and by corollary \ref{c7} is
perfect over $B$. In particular, it is compact and cocompact as an
ind-pro-perfect $B$-module. Since $\HH_{DR}(X,E^\vee \otimes E)$ is
ind-perfect it is cocompact as an ind-pro-perfect $B$-module. We thus
have that $\HHcd(X,E^\vee \otimes E)$ is also cocompact as an
ind-pro-perfect $B$-module.  Since it is pro-perfect, it must be
perfect. But this implies that $\HH_{DR}(X,E^\vee \otimes E)$ is
perfect and thus that $E$ is Fredholm by corollary \ref{cd11}.  \hfill
$\Box$

\

\noindent
The previous lemma shows that $\Vect^{\nabla}_{V_\s}(X)$ has a global
cotangent complex which is furthermore perfect.  In order to apply
Theorem~\ref{t2} it remains to prove that $\Vect^{\nabla}_{V_\s}(X)$
satisfies the three conditions $(2)$, $(5)$ and $(6)$.  These three
statements are properties of the restriction of
$\Vect^{\nabla}_{V_\s}(X)$ to underived $k$-algebras. Let us denote
this restriction by $\Vect^{\nabla}_{V_\s}(X)_0$.

We start by studying the diagonal morphism of
$\Vect^{\nabla}_{V_\s}(X)_0$ in order to check condition $(2)$ of
Theorem~\ref{t2}.

\begin{lem}\label{l7}
The diagonal morphism
$$
\mathsf{diag} : \Vect^{\nabla}_{V_\s}(X)_0 \longrightarrow
\Vect^{\nabla}_{V_\s}(X)_0 \times \Vect^{\nabla}_{V_\s}(X)_0
$$
is representable by a scheme of finite type over $k$. 
\end{lem}

\

\noindent
\textit{Proof of the lemma:} The statement of the lemma is equivalent
to the statement that for any discrete cdga $B$ and any two points $E$
and $F$ in $\Vect^{\nabla}_{V_\s}(X)(B)$, the sheaf of isomorphisms
$Iso(E,F)$ is representable by a scheme of finite type over $Spec\,
B$. This sheaf is an open sub-sheaf inside the sheaf of morphisms
$Hom(E,F)$ from $E$ to $F$, it is therefore enough to prove that
$Hom(E,F)$ is representable by a scheme of finite type over $B$. The
value of this sheaf over a $B$-algebra $B'$ is given as the fiber at
the identity of the restriction map
$$
\xymatrix{ 0 \ar[r] & Hom(E,F)(B') \ar[r] &
  \HH^{0}_{DR}(X,E^\vee\otimes F \otimes_B B') \ar[r] & \HH^{0}_{DR}(\fb
  X,V_{\s}^\vee\otimes V_\s)\otimes_k B'.}
$$
In other words $Hom(E,F)$ is the sheaf of morphisms with compact
supports (i.e. restrict to the identity morphism on $\fb X$) from $E$
to $F$. Because $E^\vee\otimes F$ is automatically Fredholm, the
functor sending $B'$ to $\HH^{0}_{DR}(X,E^\vee\otimes F \otimes_B B')$ is
the $H^0$-functor of a perfect complex over $B$ of amplitude $[0,\s)$,
    and thus is representable by a scheme of finite type.

\begin{sublem}\label{sl}
Let $K$ be a perfect complex on a commutative $k$-algebra $B$, and
suppose that $K$ has amplitude contained in $[0,\s)$. Then the functor
    $B' \mapsto H^0(K\otimes_B B')$ is representable by an affine
    scheme of finite presentation over $B'$.
\end{sublem}

\

\noindent
\textit{Proof of the sublemma:} Because of the amplitude hypothesis
$K$ can be presented by a bounded complex of projective modules of
finite rank
$$
\xymatrix{
  0 \ar[r] & K^0 \ar[r] & K^{1} \ar[r] & \cdots \ar[r] & K^n}
$$
for some integer $n$. The functor under consideration is then the
kernel of $K^0 \longrightarrow K^1$, that is the kernel of a morphism
between vector bundles over $Spec\, B$, and the result follows as
affine schemes of finite presentation over $B$ are stable by fiber
products.  \hfill $\Box$

\

\noindent
Going back to the proof of lemma \ref{l7}, the sublemma and the fact
that $E^\vee\otimes F$ is automatically Fredholm, imply that the two
functors
$$B' \mapsto \HH^{0}_{DR}(X,E^\vee\otimes F \otimes_B B')$$
and
$$B' \mapsto \HH^{0}_{DR}(\fb X,V_{\s}^\vee\otimes V_\s)\otimes_k B'$$ are
representable by affine schemes of finite presentation over $Spec\,
B$.  We thus get that the sheaf $Hom(E,F)$ is also representable by an
affine scheme of finite presentation over $Spec\, B$, which completes
the proof of the lemma.  \hfill $\Box$

\

The previous lemma implies that condition $(2)$ of Theorem
\ref{t2} is also satisfied. Indeed, the diagonal morphism has the
property that it is nil-complete, inf-cartesian and possesses a
perfect cotangent complex, so the fact that it is representable on the
level of truncations implies that it is representable (see
\cite{hagII}).  The condition $(1)$ of Theorem \ref{t2} is also
satisfied as no components of $X$ are assumed to be proper, so for any
$V\in \Vect^{\nabla}(X)(k)$ the induced morphism
$$
\HH^{0}_{DR}(X,V^\vee \otimes V) \longrightarrow
\HH^{0}_{DR}(\fb X,R(V)^\vee\otimes R(V))$$
is injective. 
It thus remains to check conditions
$(5)$ and $(6)$ of Theorem \ref{t2}. 

First we will check that condition $(5)$ of Theorem \ref{t2} is
satisfied by $\Vect^{\nabla}_{V_\s}(X)$. By \cite{mo} we can chose a
(possibly stacky) good compactification $X \hookrightarrow \Xbar$ such
that the underlying bundle of $V$ extends to a vector bundle
$\barV$ on $\Xbar$.  We denote by $D \hookrightarrow \Xbar$ the
divisor at infinity.  The connection on $V$ can then be represented
by a connection with poles
$$
\barnabla : \barV \longrightarrow \Omega_{\Xbar}^1(nD)
\otimes_{\OO_{\Xbar}}\barV
$$
for some integer $n$. The morphism $\barnabla$ can also be
interpreted as a splitting of the Atiyah extension with poles along
$D$:
$$
E(\barV,n) : \xymatrix{
0 \ar[r] & \Omega_{\Xbar}^1(nD) \otimes_{\OO_{\Xbar}}\barV 
\ar[r] & \mathcal{P}(\barV)(nD) \ar[r]
& \barV \ar[r] & 0,}
$$
where $\mathcal{P}(\barV)(nD)$ is the vector bundle of
principal parts of $\barV$ possibly with poles of order at most
$n$ along $D$.

We consider the (underived) stack of pairs $(\barW,\delta)$,
consisting of a vector bundle $\barW$ on $\Xbar$ and a flat connection
$\delta$ on $\barW$ with poles of order at most $n$ along $D$.  By
definition this stack sends a commutative $k$-algebra $B$ to the
groupoid of vector bundles $\barW$ on $\Xbar\times Spec\, B$, together
with a splitting $\delta$ of the exact sequence of bundles on $\Xbar\times
Spec\, B$:
$$
E(\barW,n) : \xymatrix{
0 \ar[r] & \Omega_{\Xbar}^1(nD) \otimes_{\OO_{\Xbar}}\barW 
\ar[r] & \mathcal{P}_{\Xbar,B}(\barW)(nD)  \ar[r]
& \barW  \ar[r] & 0,}
$$
satisfying the integrability condition $\delta^2=0$ as a section of
$\Omega_{\Xbar}^2(2nD)\otimes_{\OO_{\Xbar}}End(\barW)$. Here
$\mathcal{P}_{\Xbar,B}(\barW)(nD)$ denotes the sheaf o principal parts
of $\barW$, taken relative to the map $\Xbar\times Spec\, B \to Spec\,
B$ and with poles of order at most $n$ along $D \times Spec\, B$.

Let us denote
this stack by $\mathcal{F}_{\Xbar}$. This is clearly an Artin stack
locally of finite type over $k$. In the same manner we can define
$\mathcal{F}_{X}:=\Vect^{\nabla}(X)_0$ - the underived stack of vector
bundles with flat connections on $X$, as well as
$\mathcal{F}_{\widehat{\Xbar}}$ - the stack of vector bundles
on the formal completion $\widehat{\Xbar}$ endowed with flat connections
with poles of order at most $n$ along $D \hookrightarrow
\widehat{\Xbar}$.  Finally, we have $\mathcal{F}_{\fb
  X}:=\Vect^{\nabla}(\fb X)_0$. The formal gluing of \cite{hpv} again
implies that there exists a cartesian square of underived stacks
$$
\xymatrix{
\mathcal{F}_{\Xbar} \ar[r] \ar[d] & \mathcal{F}_{X} \ar[d] \\
\mathcal{F}_{\widehat{\Xbar}} \ar[r] & \mathcal{F}_{\fb X}.}
$$
The stack $\mathcal{F}_{\widehat{\Xbar}}$ is a limit of Artin stacks
locally of finite type, and thus satisfies the condition $(5)$ of
Theorem \ref{t2}.  The stack $\mathcal{F}_{\Xbar}$ satisfies the
conditions $(5)$ and $(6)$ of theorem \ref{t2}.  This implies that the
fiber of the left vertical  map, taken at
$\widehat{\jfrak}^{*}(\barV,\barnabla)$, will satisfy condition $(5)$.  But
by construction this fiber is the truncated stack
$\Vect^{\nabla}_{V_\s}(X)_0$.  This implies that
$\Vect^{\nabla}_{V_\s}(X)_0$ satisfies the condition $(5)$ of the
theorem \ref{t2}, as desired.

\

\noindent
Finally we need to show that $\Vect^{\nabla}_{V_\s}(X)$ satisfies
condition $(6)$ of Theorem \ref{t2}. For this, let $B= \op{colim}_i
B_i$ as in $(6)$, and assume that each $B_i$ as well as $B$ are
noetherian rings. We consider
$$
\op{colim}_i \Vect^{\nabla}_{V_\s}(X)(B_i) \longrightarrow 
\Vect^{\nabla}_{V_\s}(X)(B).
$$
By Lemma \ref{l7} this map is injective and so we need to show it is
surjective as well.  Let us fix an object in
$\Vect^{\nabla}_{V_\s}(X)(B)$, represented by a pair $(E,\alpha)$, of
$E \in \Vect^{\nabla}(X)(B)$ and $\alpha : \mathsf{R}(E) \simeq V_\s \times_k
B$ in $\Vect^{\nabla}(\fb X)(B)$.  Since the stack $\Vect^{\nabla}(X)$
of flat bundles on $X$ is locally of finite presentation, there is an
$i$ and $E_i \in \Vect^{\nabla}(X)(B_i)$ such that $E_i
\otimes_{B_i}B\simeq E$.

We now consider the sheaf $\mathcal{I}$ of 
isomorphisms between $\mathsf{R}(E_i)$ and $V_\s\otimes_k B_i$, which is
a sheaf on the big \'{e}tale site of affine schemes over $S_i=Spec\, B_i$.
This sheaf is a subsheaf in $\mathcal{J}$ the sheaf of all 
morphisms from $\mathsf{R}(E_i)$ to $V_\s\otimes_k B_i$. 

\begin{lem}
There exists a non-empty Zariski open $U_i \subset S_i=Spec\, B_i$
such that the restriction of the sheaf $\mathcal{J}$ is representable
by a scheme of finite type over $U_i$.
\end{lem}

\

\noindent
\textit{Proof of the lemma:} This is similar to the argument we used
in Lemma \ref{l7}. We have to prove that if we set
$$
E_i':= \mathsf{R}(E_i)^\vee \otimes V_\s\otimes_k B_i \in
\Vect^{\nabla}(\fb X)(B_i),
$$ then the Tate object $\widetilde{\HH}_{DR}(\fb X,E_i')[f^{-1}]$ is
a perfect $B_i[f^{-1}]$-module, for some non-zero localization
$B_i[f^{-1}]$.  For this we use the criterion from Corollary \ref{c6}
and Proposition \ref{phol}.

First let us recall some notation. We will again write $j : X
\hookrightarrow \Xbar$ be the embedding in the good
compactification.  Recall also that at the beginning of the proof
of Theorem~\ref{t1} we fixed a flat bundle $V \in \Vect^{\nabla}(X)$
satisfying $\mathsf{R}(V) \simeq V_\s$.

We first notice that $j_* E$ is a perfect $\D_{\Xbar,B}$-module
on $\Xbar\times Spec\, B$. This is a local statement on $\Xbar$ which
reduces to the following algebraic fact. Let $A$ be a smooth
$k$-algebra of finite type and $f\in A$.  We consider
$\widehat{A\otimes_k B}$, the formal completion of $A\otimes_k B$ at
$f\otimes 1$. We denote by $\widehat{\D}_{\Xbar,B}$ the ring of
completed relative differential operators. As a module it is
$\widehat{A\otimes_k B}_{A\otimes_k B}(\D_{\Xbar} \otimes_k B)$, where
the ring structure is defined naturally by making $\D_{\Xbar}
\otimes_k B$ act on the completion $\widehat{A\otimes_k B}$ by
extending derivations to the completion.  In the same manner we let
$\widehat{\D}_{X,B}$ be $\widehat{\D}_{\Xbar,B}[f^{-1}]$. Using the
formal gluing of \cite{bh} we have a cartesian square of
$\s$-categories
$$
\xymatrix{
\Dqcoh(\D_{\Xbar,B}) \ar[r]^-{j^{*}} \ar[d] & \Dqcoh(\D_{X,B}) \ar[d] \\
\Dqcoh(\widehat{\D}_{\Xbar,B}) \ar[r]_-{\Jfb^{*}}  &
\Dqcoh(\widehat{\D}_{X,B}).}
$$
The functor $\Jfb^{*}$ is an explicit incarnation of the
pullback functor for $\mathcal{D}$-modules from the full infinitesimal
neighborhood $\widehat{\Xbar}$ of $D$ in $\Xbar$ to the ``punctured
infinitesimal neighborhood'' $\fb X$ of $D$ in $\Xbar$. Similarly to
the way the map $j : X \hookrightarrow \Xbar$ gives rise to the
pullback/pushforward adjoint pair of functors $j^{*} \dashv j_{*}$
acting on relative $\mathcal{D}$-modules, the functor $\Jfb^{*}$
posseses a right adjoint denoted by $\Jfb_{*}$. This right adjoint
is given by considering a $\widehat{\D}_{X,B}$-module
as a $\widehat{\D}_{\Xbar,B}$-module via the canonical 
morphism of sheaves of dg-algebras $\widehat{\D}_{\Xbar,B} \to 
\widehat{\D}_{X,B}$.

Therefore, for $j_* E \in \Dqcoh(\D_{\Xbar,B})$ to be perfect it
is enough that its restrictions as $\D_{X,B}$ and
$\widehat{\D}_{\Xbar,B}$ modules are both perfect. But, the first of
these restrictions is $E$ which is perfect over $\D_{X,B}$, and the
second of these restrictions corresponds to
$\Jfb_{*}(V_\s^\vee \otimes_k
B)$. This is perfect because it is the restriction to $\widehat{\Xbar}$ of
$j_*(V^\vee)\otimes_k B \in \Dqcoh(\D_{\Xbar,B})$, which is perfect
because of Bernstein's theorem asserting that $j_*(V^\vee)$ is a
coherent and holonomic complex of $\D_{\Xbar}$-modules.

We thus have that $j_* E$ is a perfect
$\D_{\Xbar,B}$-module. Moreover it is also holonomic (see Proposition
\ref{phol}). Indeed, because $\mathsf{R}(E)$ is isomorphic to
$V_\s\otimes_k B$, its characteristic cycle is contained in \linebreak
$\Lambda \times Spec\, B \subset T^\vee X \times Spec\, B$, where
$\Lambda=Char(\Jfb_*(V_\s))$. Since the $\s$-functor sending $B$ to perfect
$\D_{\Xbar,B}$-modules is locally of finite presentation, we can chose
$i$ and $F_i \in \D_{\mathsf{perf}}(\D_{\Xbar,B})$ so that
$j_* E\simeq F_i \otimes_{B_i}B$. By enlarging $i$ if necessary,
we can also assume that the characteristic variety of $F_i$ is
contained in $\Lambda \times Spec\, B_i$, and thus that $F_i$ is
moreover holonomic. Also, we can assume that $F_i$ and $j^*E_i^\vee$
are isomorphic as objects in $\Vect^{\nabla}(X)(B_i)$.

As now both $F_i$ and $\Jfb_*(V_\s)\otimes_k B_i$ are perfect and
holonomic, Proposition $\ref{phol}$ implies that \linebreak
$F_i\otimes_{\OO}(\Jfb_*(V_\s^\vee)\otimes_k B_i)\simeq j_*(E_i')$ remains
perfect over $\D_{\Xbar,B_i[f^{-1}]}$ for some non-zero localization
of $B_i$.  Working with $V_\s^\vee$ and $E^\vee$ from the start we
prove the same manner that $j_*((E_i')^\vee)$ is also perfect. By
Corollary \ref{c6} this implies that $\widetilde{\HH}_{DR}(\fb
X,E_i')$ is perfect which finishes the proof of the lemma.  \hfill
$\Box$

\

\noindent
By the above lemma $\mathcal{J}$ is representable by a scheme of
finite type.  The sheaf $\mathcal{I}$ clearly is an open subsheaf of
$\mathcal{J}$ and thus is also representable by a scheme of finite
type over an non-empty open on $Spec\, B_i$.  The canonical
isomorphism $\alpha : \mathsf{R}(E) \simeq V_\s \otimes_k B$, which an element
in $\mathcal{I}(B)$ is then definable over some $B_i[f^{-1}]$ for some
$i$ and non-zero localization, say $\alpha_i : \mathsf{R}(E_0)[f^{-1}] \simeq
V_\s \otimes_k B_i[f^{-1}]$.  The pair $(E_0,\alpha_i)$ defines an
object in $\Vect^\nabla(X)_{V_\s}(B_i[f^{-1}])$ whose image in
$\Vect^{\nabla}_{V_\s}(X)(B[f^{-1}])$ is the restriction of our
original object $E$.

This finishes the proof of condition
$(6)$ of theorem \ref{t2}, and thus of theorem \ref{t1}.
\hfill $\Box$

\

\noindent
Unfortunately, we do not know if Theorem \ref{t1} can be
strengthened to the statement that $\Vect^{\nabla}_{V_\s}(X)$ is
representable by an algebraic space locally of finite type over
$k$. The only missing condition would be that
$\Vect^{\nabla}_{V_\s}(X)$ is also locally of finite presentation, a
condition that we havent been able to prove or disprove.

\newpage
\appendix

\Appendix{Derived quasi-algebraic spaces and Artin's representability}
\label{sec:artin-lurie}

\setcounter{equation}{0}

\

\noindent
In this section we have gathered some definitions and results on
\emph{\bfseries derived quasi-algebraic spaces} and the corresponding
representability criterion. Derived quasi-algebraic spaces are slight
gneralizations of derived algebraic spaces for which atlases only
exist generically. These derived stacks are not algebraic in general,
but are algebraic as soon as the functors they represent are locally
of finite presentation.

To make sense of such spaces, we will need the following notion of a
dominant morphism to a not necessarily algebraic derived stack
$F$. Assume that $F$ is a derived stack which has a perfect global
cotangent complex, and is nil-complete and infinitesimally
cartesian. We will also assume that $F$ is integrable, that is for any
local complete noetherian discrete $k$-algebra $A= \lim_i
A/\mathfrak{m}^i$, the natural morphism
$$F(A) \longrightarrow \lim_i F(A/\mathfrak{m}^i)$$
is bijective.

For any such $F$, any field $K$ which is finitely generated over
$k$, and any point $x : Spec\, K \longrightarrow F$, there exists by
\cite[Theorem~18.2.5.1]{lu2} a complete local noetherian cdga $A$ with
residue field $K$, and a formally smooth morphism
$$Spf(A) \longrightarrow F$$
extending the point $x$. We get this way a morphism from its truncation
$$Spf(\pi_0(A)) \longrightarrow F,$$
and by integrability a well defined morphism
$$\widehat{x} : Spec(\pi_0(A)) \longrightarrow F.$$
A morphism $\widehat{x}$ obtained this way will be called
a \emph{\bfseries formally smooth lift of $x$}. 

\begin{df}\label{ddom}
For a derived stack $F$ as above and a derived scheme $X$ locally of
finite presentation over $k$, with a morphism $f : X \longrightarrow
F$. We say that $f$ is \emph{\bfseries dominant} if for any finitely
generated $k$-field $K$, any point $x : Spec\, K \longrightarrow F$,
and any formally smooth lift $\widehat{x} : Spec(\pi_0(A)) \longrightarrow
F$, the derived scheme $X \times_{F} Spec(\pi_0(A))$ is non-empty.
\end{df}

\

\noindent
Note that if $F$ is itself representable by a derived algebraic space
locally of finite presentation, then $f : X \longrightarrow F$ is
dominant in the sense above if and only if for any \'{e}tale morphism
$Spec\, B \longrightarrow F$ we have $X\times_F Spec\, B \neq
\varnothing$.  Indeed, assume that there is an \'{e}tale map $Spec\, B
\longrightarrow F$ whose pull-back to $X$ is empty.  We pick a point
$x$ of $Spec\, B$ and consider the corresponding  formal completion
$\widehat{B}_{x}$ of $B$.  Since $Spec\, B \longrightarrow F$ is \'{e}tale
the composition
$$
Spec\, \pi_0(\widehat{B}_{x}) \longrightarrow Spec\, B \longrightarrow F
$$
is a formally smooth lift of $x$. By construction the pull-back
$Spec\, \pi_0(\widehat{B}_{x}) \times_F X$ is empty. This shows that
the above notion of dominant map is a generalization of the notion of
a morphism with Zariski dense image.

We can now give the definition of a derived quasi-algebraic spaces 
as derived stacks with dominant smooth atlases as follows.

\begin{df}\label{d17}
A derived stack $F$ is a \emph{\bfseries derived quasi-algebraic
  space} (locally of presentation with schematic diagonal of finite
presentation) if it satisfies the following conditions.

\begin{enumerate}
\item[(i)] The diagonal of the stack $F \longrightarrow F \times F$ is
  representable by a derived scheme of finite presentation.
\item[(ii)] The derived stack $F$ has a perfect global cotangent
  complex, and is nil complete and infinitesimally cartesian.
\item[(iii)] The derived stack $F$ is integrable: for any local
  complete noetherian discrete $k$-algebra $A=\lim_i
  A/\mathfrak{m}^i$, the natural morphism
$$F(A) \longrightarrow \lim_i F(A/\mathfrak{m}^i)$$
is bijective. 
\item[(iv)] There exists a family of cdga $A_i$ of finite presentation
  over $k$ and a morphism $p : \sqcup Spec\, A_i \longrightarrow F$
  such that
\begin{enumerate}
\item[(a)] For each $i$ the morphism $Spec\, A_i \longrightarrow F$ is
  smooth.
\item[(b)] The morphism $p$ is dominant in the sense of definition
  \ref{ddom} above.
\end{enumerate} 
\end{enumerate}
\end{df}

\

\noindent
A derived quasi-algebraic space is algebraic if and only if the
functor $F$ is futhermore locally of finite presentation. This follows
from Artin-Lurie's representability theorem
\cite[Theorem~18.3.0.1]{lu2}. Similarly derived quasi-algebraic spaces
can be characterized by the following version of Artin's
representability.

\begin{thm}\label{t2}
A derived stack $F$ is a derived quasi-algebraic space 
if it satisfies the following conditions.

\begin{enumerate}

\item[(1)] For any discrete cdga $B$ the simplicial set $F(B)$ is
  $0$-truncated.
\item[(2)] The diagonal morphism of its truncation is representable by
  a scheme of finite presentation.
\item[(3)] The derived stack $F$ has a perfect global cotangent complex.
\item[(4)] The derived stack $F$ is nil-complete and infinitesimally
  cartesian.
\item[(5)] For any discrete local $k$-algebra $(A,\mathfrak{m})$
  essentially of finite type, with completion $\widehat{A}=\lim_i
  A/\mathfrak{m}^i$, the morphism $F(\widehat{A}) \longrightarrow
  \lim_i F(A/\mathfrak{m}^i)$ is an equivalence.
\item[(6)] For any filtered system of noetherian discrete commutative
  $k$-algebras $B=\op{colim}_i B_i$ and any $x\in F(B)$, there exists
  an index $i$ and a non-empty Zariski open $U_i \subset Spec\, B_i$
  with $U=U_i\times_{Spec\, B_i}Spec\, B$ non-empty, and such that the
  restriction of $x$ lies in the image of $F(U_i) \longrightarrow
  F(U)$.
\end{enumerate}
\end{thm}

\

\noindent
\textit{Sketch of a proof:} The proof is essentially the same as the
usual representability theorem in \cite{lu2}.

Consider fields $K$ which are finitely generated over $k$. 
For any morphism 
$$x : Spec\, K \longrightarrow F$$
we can use \cite[Theorem~18.2.5.1]{lu2}
to find a local complete and noetherian cdga $(A,\mathfrak{m})$ 
with residue field $K=A/\mathfrak{m}$ and a factorisation
$$Spec\, K \hookrightarrow Spf(\widehat{A}) \longrightarrow F,$$
where the second map is formally smooth (i.e. its 
relative contangent complex is a vector bundle). We write 
$B=\pi_0(\widehat{A})$, which is a complete local discrete $k$-algebra 
with residue field $K$, and 
consider
the induced morphism on the 
truncation $\widehat{x} : Spf(B) \longrightarrow F$.
We can use condition (4) to lift this to a factorization
$$Spec\, K \hookrightarrow Spec(B) \longrightarrow F.$$

As explained in the proof of \cite[Theorem~18.2.5.1]{lu2}, there
exists a $k$-algebra of finite type $B' \subset B$, such that if
$p=\mathfrak{m}\cap B'$, then the induced morphism on formal
completions
$$\widehat{B}'_p \longrightarrow B$$ is surjective (take $A'$ big
enough so that it contains generators for $K$ over $k$ as well as
generators of the $k$-vector space $\mathfrak{m}/\mathfrak{m}^2$).  We
can now apply Popescu's theorem to the regular morphism $B'
\longrightarrow \widehat{B}'_p$ and thus write $\widehat{B}'_p=
\op{colim}_i B'_i$ as a filetered colimit of smooth $B'$-algebras. Since
$B$ is finitely presented as a $\widehat{B}'_p$-algebra, we can find an
index $i$ and a $B'_i$-algebra $C'_i$ of finite presentation such that
$$C\simeq \op{colim}_i (\widehat{B}'_{p}B\otimes_{B'_i}C'_i).$$ We let
$C_i:=\widehat{B}'_{p}B\otimes_{B'_i}C'_i$, which is a $B'_i$-algebra
of finite presentation, and thus is itself of finite presentation over
$k$.

We now apply condition $(6)$ to the morphism $Spec\, B \longrightarrow
F$, and get that there exists an integer $i$ and a Zariski open
$U_i=Spec\, C_i[f^{-1}] \subset Spec\, C_i$, with $U=Spec B[f^{-1}]$
non-empty, and which fits in a commutative diagram
$$
\xymatrix{
U_i  \ar[dr] & \\
U \ar[u] \ar[d] & F \\
Spec\, B. \ar[ru] & }
$$

\begin{lem}
With the notation above, and enlarging $i$ is necessary, the morphism
$p : U_i \longrightarrow F$ constructed above is formally smooth is
the underived sense: $\tau_{\leq -1}(\mathbb{L}_{U_i/F})$ is a vector
bundle in degree $0$.
\end{lem}

\

\noindent
\textit{Proof of the lemma:} First of all $U_i$ being of finite type
together with the fact that the diagonal of $F$ is representable of
locally of finite presentation implies that $p$ is representable and
locally of finite type in the underived sense. It thus only remains to
show that $p$ is also formally smooth in the underived sense,
i.e. that its relative $1$-truncated cotangent complex $\tau_{\leq
  1}(\mathbb{L}_{U_i/F})$ is a vector bundle.

For this we first notice that $\mathbb{L}_{U_i/F}$ is almost perfect
(i.e. quasi-isomorphic to a complex of free modules of finite rank
over $C_i[f^{-1}]$ concentrated in degree $(-\s,0]$). Since we are
only interested in its truncation $\tau_{\leq
  1}(\mathbb{L}_{U_i/F})$ we will be able to act as if
$\mathbb{L}_{U_i/F}$ is in fact perfect (simply replace it by a
perfect complex having the same cohomology in degree $[-n,0]$ for n
big enough).  We start by computing the pull-back of
$\mathbb{L}_{U_i}$ to $U=Spec\, B[f^{-1}]$.

Consider the exact triangle of complexes of $B$-modules (where
$\mathbb{L}_{A}$ stands for $\mathbb{L}_{A/k}$ for any $k$-algebra
$A$).
$$
\xymatrix{
  \mathbb{L}_{C_i}\otimes_{C_i}B \ar[r] &
  \mathbb{L}_B \ar[r] & \mathbb{L}_{B/C_i}.}
$$
Since $B$ is complete with respect to its maximal ideal
$\mathfrak{m}$, for any connective dg-module $E$ over $B$, we have its
completion $\widehat{E}:= \lim_{i} E\otimes_{B}B/\mathfrak{m}^i$,
which is another connective$B$-dg-module together with a natural
morphism $E \longrightarrow \widehat{E}$. Moreover, when $E$ is almost
perfect this morphism is a quasi-isomorphism.  We can then complete
the terms in the above triangle to get a new triangle
$$
\xymatrix{ \widehat{\mathbb{L}_{C_i}\otimes_{C_i}B} \ar[r] &
  \widehat{\mathbb{L}_B} \ar[r] & \widehat{\mathbb{L}_{B/C_i}}.}
$$
As $\mathbb{L}_{C_i}$ is almost perfect the first term is simply
$\mathbb{L}_{C_i}\otimes_{C_i}B$. Moreover, by base change
$\mathbb{L}_{B/C_i}$ is naturally equivalent to
$\mathbb{L}_{\widehat{B}'/C_i'}\otimes_{B'}B$. Since $B' \longrightarrow
B$ is a surjective local morphism we see that the base change of
$\widehat{\mathbb{L}_{\widehat{B}'/C_i'}}$, considered as a pro-object in
connective $B'$-dg-modules, by the map $\widehat{B}'\longrightarrow
B$, is the pro-object $\widehat{\mathbb{L}_{B/C_i}}$.

We now use that $\widehat{B}'$ is the completion of $B'$ along the
maximal ideal $\mathfrak{m} \subset B'$, and so for all $i$ we
have $\mathbb{L}_{\widehat{B}'/B'}\otimes_{B'}B'/\mathfrak{m}^i \simeq 0$.
We thus
have an equivalence of pro-objects
$$\widehat{\mathbb{L}_{\widehat{B}'/C_i}}\simeq
\widehat{\mathbb{L}_{C_i/B'} \otimes_B' \widehat{B}'}[1].$$
Since
$C_i$ is smooth over $B'$ we therefore conclude that the pro-object
$\widehat{\mathbb{L}_{B/C_i}}$ is a vector bundle in degree $1$, and
so its realization as a connective $B$-dg-module is
$\mathbb{L}_{C_i/B'}\otimes_{C_i}B[1]$.
Our original triangle can
therefore be written as
$$
\xymatrix{ \mathbb{L}_{C_i}\otimes_{C_i}B \ar[r] &
  \widehat{\mathbb{L}_B} \ar[r] & V[1]}
$$
with $V$ a vector bundle on $Spec\, B$. We can now localize this
triangle to the open $U=Spec\, B[f^{-1}]$ in order to get new triangle
on $U$
$$
\xymatrix{
  \mathbb{L}_{C_i[f^{-1}]}\otimes_{C_i[f^{-1}]}B[f^{-1}] \ar[r] &
  \widehat{\mathbb{L}_B}[f^{-1}] \ar[r] & V[f^{-1}][1].}
$$
The morphisms $q : U\rightarrow U_i \longrightarrow F$ induces a
morphism
$$\xymatrix{
  q^*(\mathbb{L}_F) \ar[r] \ar[rd] & \mathbb{L}_{C_i[f^{-1}]}
  \otimes_{C_i[f^{-1}]}B[f^{-1}] \ar[d] \\
 & \mathbb{L}_{B}[f^{-1}],}$$
which factors throught completions since 
$q^*(\mathbb{L}_F)$ is a perfect complex by our condition $(3)$. We get 
$$\xymatrix{
  q^*(\mathbb{L}_F) \ar[r] \ar[rd] & \mathbb{L}_{C_i[f^{-1}]}
  \otimes_{C_i[f^{-1}]}B[f^{-1}]\ar[d] \\
 & \widehat{\mathbb{L}_B}[f^{-1}],}$$
and the induced morphism on the cones sits in an exact triangle
$$
\xymatrix{ \mathbb{L}_{U_i/F}\otimes_{C_{i}[f^{-1}]}B[f^{-1}] \ar[r]
  & \widehat{\mathbb{L}_{Spec\ B/F}}[f^{-1}] \ar[r] & V[1]}.
$$
Because $Spf \longrightarrow F$ was chosen to be  formally smooth in the
underived sense we have that \linebreak $\tau_{\leq
  -1}(\widehat{\mathbb{L}_{Spec\ B/F}}[f^{-1}])$ is a vector bundle in
degree $0$. The conclusion is that $\mathbb{L}_{U_i/F}$ is an almost
perfect complex over $C_i[f^{-1}]$ such that its base change to
$B[f^{-1}]=colim (C_i[f^{-1}])$ has vanishing $H^-1$ and a vector
bundle as $H^0$. This implies that the same is true for
$\mathbb{L}_{U_i/F}\otimes_{C_i[f^{-1}]}C_j[f^{-1}]$ for some big
enough $j$. \hfill $\Box$

\

\noindent
Going back to the proof of the theorem, we use again
\cite[Theorem~18.2.5.1]{lu2} but this time for $U_i \longrightarrow F$,
which by the lemma can be chosen to be formally smooth in the
underived sense. We can thus produce a smooth morphism
$$W_i \longrightarrow F$$ where $W_i$ is a derived affine scheme whose
truncations coincides with the given map $U_i \longrightarrow F$. The
derived scheme $W_i$ is itself of finite presentation over $k$ as its
truncation is of finite type and its cotangent complex is perfect
(because its maps smoothly to $F$).

Taking the union of all morphisms $W_i \longrightarrow F$ constructed
above provides the required generic atlas for $F$ as in definition
\ref{d17}. \ \hfill $\Box$

\smallskip

\noindent
Tony Pantev, {\sc University of Pennsylvania, Department of Mathematics, DRL 209 South 
33rd Street,
Philadelphia, PA 19104-6395.}\\
email:tpantev@math.upenn.edu

\bigskip

\noindent
Bertrand To\"{e}n, {\sc  Universit\'e de Toulouse  \& CNRS, IMT, 118 route de Narbonne, 31062 Toulouse, France.} \\
email:bertrand.toen@math.univ-toulouse.fr


\begin{thebibliography}{99}

\bibitem[Ba-Ca-Fi] {bacafi} 
Baldassarri, Francesco; Cailotto, Maurizio; Fiorot, Luisa 
\emph{Poincar\'e duality for algebraic de Rham cohomology.} 
Manuscripta Math. \textbf{114} (2004), no. 1, 61–116. 

\bibitem[Be-Te]{bete} Ben-Bassat, Oren; Temkin, Michael \emph{Berkovich 
spaces and tubular descent.}
Adv. Math., \textbf{234}:217--238, 2013.


\bibitem[Bh]{bh} Bhatt, Bhargav \emph{Formal glueing of module categories.}
Available at 
http://www-personal.umich.edu/$~$bhattb/math/formalglueing.pdf

\bibitem[Ca]{ca} Calaque, Damien
\emph{Lagrangian structures on mapping stacks and semi-classical TFTs.} 
in Stacks and categories in geometry, topology, and algebra, 1--23,
Contemp. Math., \textbf{643}, Amer. Math. Soc., Providence, RI, 2015. 

\bibitem[CPTVV]{cptvv} Calaque, Damien; Pantev, Tony; To\"en, 
Bertrand; Vaqui\'e, Michel; Vezzosi, Gabriele 
\emph{Shifted Poisson structures and deformation quantization.} 
J. Topol. \textbf{10} (2017), no. 2, 483--584.



\bibitem[Ef]{ef} Efimov, Alexander \emph{Categorical formal punctured
    neighborhood of infinity, I.}  Preprint arXiv:1711.00756

\bibitem[Fo-Ro]{fock-rosly} Fock, V. V.; Rosly, A. A. \emph{Poisson
    structure on moduli of flat connections on Riemann surfaces and
    the r-matrix}.  Moscow Seminar in Mathematical Physics, 67–86,
  Amer. Math. Soc. Transl. Ser. 2, 191, Adv. Math. Sci., 43,
  Amer. Math. Soc., Providence, RI, 1999.


  
\bibitem[Ga-Ro]{gr} Gaitsgory, Dennis; Rozenblyum, Nick \emph{Crystals
    and D-modules.} Pure Appl. Math. Q. \textbf{10}, No. 1, 57-154 (2014). 




\bibitem[Gol]{goldman} Goldman, William \emph{Mapping class group
    dynamics on surface group representations}. Problems on mapping
  class groups and related topics, 189–214, Proc. Sympos. Pure Math.,
  74, Amer. Math. Soc., Providence, RI, 2006.

\bibitem[GHJW]{ghjw} Guruprasad, K.; Huebschmann, J.; Jeffrey, L.;
  Weinstein, A.  \emph{Group systems, groupoids, and moduli spaces of
    parabolic bundles}.  Duke Math. J. 89 (1997), no. 2, 377–412.
  
\bibitem[Gu-Ra]{guruprasad-rajan} Guruprasad, K.; Rajan, C. S.
  \emph{Group cohomology and the symplectic structure on the moduli
    space of representations}.  Duke Math. J. 91 (1998), no. 1,
  137–149.

\bibitem[He]{he} Hennion, Benjamin \emph{Tate objects in stable
  ($\s$,1)-categories.} Homology Homotopy Appl. \textbf{19} (2017),
  no. 2, 373-395.

\bibitem[He-Po-Ve]{hpv} Hennion, Benjamin; Porta, Mauro; Vezzosi, Gabriele
\emph{Formal gluing along non-linear flags.} Preprint arXiv:1607.04503

\bibitem[Lu1]{lu} Lurie, Jacob \emph{Quasi-Coherent Sheaves and
  Tannaka Duality Theorems.}  Available at
  http://www.math.harvard.edu/$~$lurie/papers/DAG-VIII.pdf

\bibitem[Lu2]{lu2} Lurie, Jacob \emph{Spectral Algebraic Geometry.}
Available at http://www.math.harvard.edu/$~$lurie/papers/SAG-rootfile.pdf

\bibitem[Lu3]{lu3} Lurie, Jacob \emph{Formal moduli problems.} Available at 
http://www.math.harvard.edu/$~$lurie/papers/DAG-X.pdf


\bibitem[Me-Sa]{mesa} Melani, Valerio; Safronov, Pavel \emph{Derived
  coisotropic structures II}. Sel. Math., New Ser. \textbf{24}, No. 4, 3119-3173 (2018). 

\bibitem[Mo]{mo} Mochizuki, Takuro \emph{Good formal structure for meromorphic flat 
connections on smooth projective surfaces.}
Algebraic analysis and around, 223–253,
Adv. Stud. Pure Math., 54, Math. Soc. Japan, Tokyo, 2009.

\bibitem[Ni]{ni} Nitsure, Nitin \emph{
Moduli of regular holonomic  D-modules with normal crossing singularities.} 
Duke Math. J. 99, No. 1, 1-39 (1999). 



\bibitem[Pa-To]{ptbetti} Pantev, Tony; To\"en, Bertrand \emph{Poisson geometry  of the  
moduli of local systems on smooth varieties.} To appear in 
Publications of RIMS, Volume \textbf{57}, issue 3-4 (2021)

\bibitem[PTVV]{ptvv}  Pantev, Tony; 
To\"en, Bertrand; Vaqui\'e, Michel; Vezzosi, Gabriele 
\emph{Shifted symplectic structures.}
Publ. Math. Inst. Hautes Etudes Sci. \textbf{117} (2013), 271-328.

\bibitem[Ra]{ra} Raskin, Sam \emph{On the notion of spectral decomposition in local 
geometric Langlands.} arXiv:1511.01378 



\bibitem[To1]{azto} To\"en, Bertrand \emph{Derived Azumaya algebras and genertors for 
twisted derived categories.} Invent. Math. \textbf{189} (2012), no. 3, 581–652.

\bibitem[To2]{to3}  To\"en, Bertrand \emph{Structures symplectiques et de Poisson sur les 
champs en cat\'egories.} arXiv:1804.10444  

\bibitem[To-Va]{tv} To\"en, Bertrand; Vaqui\'e, Michel \emph{Moduli of objects in dg-categories.}
Ann. Sci. de l’ENS Volume \textbf{40} (2007) Issue 3, Pages 387-444. 

\bibitem[To-Ve]{hagII} To\"en, Bertrand; Vezzosi, Gabriele 
\emph{Homotopical Algebraic Geometry II: geometric stacks and applications.} 
Mem. Amer. Math. Soc. \textbf{193} (2008), no. 902, x+224 pp.

\end{thebibliography}
\end{document}